\documentclass[20pt,a4paper]{article}
\usepackage{color,amsmath,amssymb, lscape,graphicx,ifthen,pdfpages}
\DeclareSymbolFont{boldsymbols}{OMS}{cmsy}{b}{n} 
\DeclareSymbolFontAlphabet{\mathbfcal}{boldsymbols} 
\usepackage[colorlinks=true,
linkcolor=webgreen,
filecolor=webbrown,
citecolor=webgreen]{hyperref}
\definecolor{bittersweet}{rgb}{1.0, 0.44, 0.37}
\definecolor{webgreen}{rgb}{0,.5,0}\definecolor{webbrown}{rgb}{.6,0,0}
\definecolor{antiquebrass}{rgb}{0.8, 0.58, 0.46}
\definecolor{antiquefuchsia}{rgb}{0.57, 0.36, 0.51}
\definecolor{apricot}{rgb}{0.98, 0.81, 0.69}
\definecolor{aquamarine}{rgb}{0.5, 1.0, 0.83}
\definecolor{atomictangerine}{rgb}{1.0, 0.6, 0.4}
\definecolor{auburn}{rgb}{0.43, 0.21, 0.1}
\definecolor{azure(colorwheel)}{rgb}{0.0, 0.5, 1.0}
\definecolor{battleshipgrey}{rgb}{0.52, 0.52, 0.51}
\definecolor{beaublue}{rgb}{0.74, 0.83, 0.9}
\definecolor{bleudefrance}{rgb}{0.19, 0.55, 0.91}

\newcommand{\seqnum}[1]{\href{http://oeis.org/#1}{\underline{#1}}}

\def\dstyle#1{$\displaystyle #1 $}
\def\noin{\noindent}
\def\pn{\par\noindent}
\def\ps{\par\smallskip}
\def\psn{\par\smallskip\noindent}
\def\pbn{\par\bigskip\noindent}
\def\Beq{\begin{equation}}
\def\Eeq{\end{equation}}
\def\Beqarray{\begin{eqnarray}}
\def\Eeqarray{\end{eqnarray}}
\def\Eq#1{eq.\,$(#1)$}
\def\sspgeq{\,\geq} 
\def\sspleq{\, \leq \,}
\def\sspkl{\, < \,}
\def\sspgr{\, > \,}
\def\sspeq{\, =\,}
\def\speq{\ =\ }
\def\sspdef{\, :=\,}
\def\spdef{\ :=\ }
\def\sspfed{\, =:\,}
\def\sspin{\, \in \,}

\def\sspcdot{\,\cdot\,}
\def\sspcirc{\,\circ\,}

\def\sspp{\, +\ }
\def\sspm{\, -\ }

\def\sspto{\,\to\,}

\def\sspcdot{\,\cdot\,}
\def\sspneq{\, \neq \,}
\def\sspequiv{\,\equiv\,}

\def\binomial#1#2{{#1} \choose {#2}}
\def\ogf{{\it o.g.f.\ }}

\def\egf{{\it e.g.f.\ }}

\def\ie{{\it i.e.},\, }
\def\eg{{\it e.g.},\, }

\def\viz{{\it viz}\, }
\def\via{{\it via}}

\def\lhs{{\it l.h.s.\, }}
\def\rhs{{\it r.h.s.\, }}

\def\op#1{{\bf #1}}
\def\floor#1{\left\lfloor{#1}\right\rfloor}

\def\spfed{\ =:\ }
\def\sspfed{\, =:\, }

\def\Simn{\,{\lower2pt\hbox{$\buildrel {\lower3pt\hbox{$n$}} \over \sim$}}\,}
\def\simn1#1{\,{\lower2pt\hbox{$\buildrel {\lower3pt\hbox{$#1$}} \over \sim$}}\,}

\def\Chi{\raise2.5pt\hbox{$\chi$}} 
\def\ssptimes{\,\times\,}
\normalbaselines
\begin{document}
\bibliographystyle{unsrt}
\rightline{Karlsruhe} \par\smallskip\noindent
\rightline{July 13 2017}
\vbox {\vspace{6mm}}
\begin{center}
{\Large {\bf On Sums of Powers of Arithmetic Progressions, and Generalized Stirling, Eulerian and Bernoulli numbers}}\\ [9mm]
Wolfdieter L a n g \footnote{ 
\href{mailto:wolfdieter.lang@partner.kit.edu}{\tt wolfdieter.lang@partner.kit.edu},\quad 
\url{http://www.itp.kit.edu/~wl}
                               } \\[3mm]
\end{center}
\vspace{2mm}
\begin{abstract}
\par\smallskip\noindent
For finite sums of non-negative powers of arithmetic progressions the generating functions (ordinary and exponential ones) for given powers are computed. This leads to a two parameter generalization of {\sl Stirling} and {\sl Eulerian} numbers. A direct generalization of Bernoulli numbers and their polynomials follows. On the way to find the {\sl Faulhaber} formula for these sums of powers in terms of generalized {\sl Bernoulli} polynomials one is led to a one parameter generalization of {\sl Bernoulli} numbers and their polynomials. Generalized {\sl Lah} numbers are also considered.     
\end{abstract}
\section{Introduction and Summary}
{\bf A) Generating functions of power sums and powers. Generalized Stirling2 and Eulerian numbers}.
\ps   
Finite sums of non-negative powers of positive integers have been studied by many authors. See {\sl Edwards}, \cite {Edwards1}, \cite {Edwards2} and {\sl Knuth} \cite{Knuth} for some history, and the books on {\sl Johannes Faulhaber} by {\sl Hawlitschek} \cite{Hawlitschek} and {\sl Schneider} \cite{Schneider}. \psn
We are interested in finite sums of power of arithmetic progressions ($PS$ for power sums)
\Beq
{\fbox{\color{blue} $PS(d,a;n,m)$}}\sspdef \sum_{j=0}^m\, (a\sspp d\,j)^n\, {\text with}\ \ n\sspin \mathbb N_0,\ m\sspin \mathbb N_0,\  d \sspin \mathbb N,\  a\sspin N_0\, .  
\Eeq  
Note that the lower summation index for $j$ is $0$. We put $0^0\sspdef 1$ if $a\sspeq 0$ and $n\sspeq 0$. (In Maple 13 \cite{Maple}  $0^0$ is put to $0$).    
It is sufficient to consider $a\sspeq 0$ if $d\sspeq 1$, and $a\sspin RRS(d)$ for $d\sspgeq 2$, where $RRS(d)$ denotes the smallest positive restricted residue system modulo $d$, \ie $RRS(d)\sspdef \{k\sspin RS(d)\, |\, \gcd(k,d) \sspeq 1\}$ with $RS(d)\sspdef\{0,\,1,\,...,\,d-1\}$, the smallest non-negative residue system modulo $d$. \psn
The aim of the first part of this paper is is to compute the ordinary (\ogf, symbolized by $G$) and exponential generating functions (\egf, symbolized by $E$) for given powers $n$. Such functions are considered in the framework of formal power series, without considering questions of convergence. Proofs will be given in section 2.
\psn
The \ogf (indeterminate $x$) is
\Beq
{\fbox{{\color{blue}$GPS(d,a;n,x)$}}}\sspdef \sum_{m=0}^{\infty}\,PS(d,a;n,m)\, x^m,\ \ n\sspin \mathbb N_0\ . 
\Eeq
The \egf (indeterminate $t$) is
\Beq
{\fbox{{\color{magenta}$EPS(d,a;n,t)$}}}\sspdef \sum_{m=0}^{\infty}\,PS(d,a;n,m)\, \frac{t^m}{m!},\ \ n\sspin \mathbb N_0\ . 
\Eeq 
As is known, the \egf is obtained from  the \ogf \via\ inverse {\sl Laplace} transform as
\Beq
E(t)\sspeq {\cal L}^{-1}\left[ \frac{1}{p}\,G\left(\frac{1}{p}\right)\right]\ ,
\Eeq 
and {\it vice versa} by a direct {\sl Laplace} transform to get $G$ from $E$.\psn
Of course, application of the binomial theorem immediately leads, after an exchange of the two finite sums, to a formula for $PS(d,a;n,m)$ in terms of the ordinary power sums $PS(n,m) \sspeq PS(1,0;n,m)$, \viz
\Beq  
PS(d,a;n,m) \sspeq \sum_{k=0}^n\, {\binomial{n}{k}}\,a^{n-k}\,d^k\,PS(k,m)\,,
\Eeq
and therefore, if we interchange an infinite sum with a finite sum,
\Beq
GPS(d,a;n,x)\sspdef \sum_{k=0}^n\, {\binomial{n}{k}}\,a^{n-k}\,d^k\,GPS(k,x)\,,  
\Eeq
with $GPS(k,x) = GPS(1,0;k,x)$.
Similarly,
\Beq
EPS(d,a;n,t)\sspdef \sum_{k=0}^n\, {\binomial{n}{k}}\,a^{n-k}\,d^k\,EPS(k,t)\,,  
\Eeq
with $EPS(k,t) = EPS(1,0;k,t)$.
Therefore it is in principle sufficient to compute $GPS(k,x)$ and use an inverse {\sl Laplace} transform to find $EPS(k,t)$. It may however be difficult (or impossible) to give its explicit form.\psn
Instead of $GPS(n,x)$ we prefer to compute the general $GPS(d,a;n,x)$ directly. In this way we find
\Beq
{\fbox{\color{blue}$GPS(d,a;n,x)$}}\sspeq \sum_{k=0}^n\, S2(d,a;n,k)\, k!\,\frac{x^k}{(1-x)^{k+2}}\ ,
\Eeq
where the generalized {\sl Stirling} numbers of the second kind (generalized subset numbers)
enter \via\ the reordering process
\Beq
(a\,{\bf 1} \sspp d\, {\bf E}_x)^n \sspfed \sum_{m=0}^n\,S2(d,a;n,m)\,x^m\, {\bf d}_x^{\ m}\,,  
\Eeq
with the {\sl Euler} operator ${\bf E}_x\sspdef x\,{\bf d}_x$ where ${\bf d}_x$ is the differentiation operator, and $\bf 1$ is the identity operator.\psn
This definition leads to the three term recurrence relation
\Beq
S2(d,a;n,m)\sspeq d\, S2(d,a;n-1,m-1) \,\sspp\, (a\sspp d\,m)\, S2(d,a;n-1,m), \ \ {\rm for}\ n\sspgeq 1\,,\ m\sspeq 0,\,1,\,...,\,n\,,  
\Eeq
with $S2(d,a;n,-1)\sspeq 0$, $ S2(d,a;n,m)\sspeq 0$ for $n\sspkl m$ and $S2(d,a;0,0)\sspeq 1$.\psn
This recurrence is obeyed by
\Beq
S2(d,a;n,m)\sspeq {\frac{1}{m!}}\, \sum_{k=0}^m\,(-1)^{m-k}\, {\binomial{m}{k}}\,(a\sspp d\,k)^n\,.  
\Eeq
These generalized {\sl Stirling} numbers build a lower triangular infinite dimensional matrix, named ${\bf S2}[d,a]$ which turns out to be an exponential convolution array like the ordinary {\sl Stirling} ${\bf S2}$ matrix, \ie a  {\sl Sheffer} matrix, denoted by 
\Beq
{\bf S2}[d,a] \sspeq (e^{a\,x},\, e^{d\,x}\sspm 1)\,.  
\Eeq
For {\sl Sheffer} matrices see the W. Lang link in $OEIS$ \cite{OEIS}, \seqnum{A006232} called ``Sheffer $a$- and $z$-sequences'', the second part, where also references are given. (Henceforth $A$-numbers will be given without quoting {\it OEIS} each time.)\psn 
A three parameter generalization of Stirling numbers of the second kind has been proposed in {\sl Bala} \cite{Bala} as $S_{(a,b,c)}$. The present generalization is $S2(d,a;n,m) \sspeq  d^m\,S_{(d,0,a)}$. There {\sl Sheffer} arrays are called exponential {\sl Riordan} arrays.\psn
A one parameter generalization is given in {\sl Luschny} \cite{Luschny1} called {\sl Stirling}-{\sl Frobenius} subset numbers, with the scaled version called there [SF-SS] with parameter m corresponding to ${\bf S2}[m,m-1]$. The [SF-S] triangle family coincides with {\sl Bala}'s $S_{(m,0,m-1)}$.   
\psn
The {\sl Sheffer} structure (exponential convolution polynomials) means that the \egf of the sequence of column $m$ is 
\Beq
ES2Col(d,a;t,m)\sspeq e^{a\,t}\,\frac{(e^{d\,t}\sspm1)^m}{m!}\,, \ \ m\sspin {\mathbb N}_0\, .
\Eeq
This corresponds to the \ogf
\Beq
GS2Col(d,a;x,m)\sspeq \frac{(d\,x)^m}{\prod_{j=0}^m\, (1\sspm (a+d\,j)\,x)}\,, \ \ m\sspin {\mathbb N}_0\ .
\Eeq
This means that the column scaled {\sl Sheffer} triangle $\widehat{S2}[d,a]\sspeq \left(e^{a\,x}, \, \frac{1}{d}\,(e^{d\,x} \sspm 1)\right)$ with elements $ \widehat{S2}(d,a;n,m)\sspeq S2(d,a;n,m)\, \frac{1}{d^m}$ are 
\Beq
\widehat{S2}(d,a;n,m)\sspeq h^{(m+1)}_{n-m}[d,a], 
\Eeq
where $h^{(m+1)}_k[d,a]$ are the complete homogeneous symmetric functions of degree $k$ of the $m+1$ symbols $a_j\sspeq a\sspp d\,j$, $j\sspeq 0,\,1,\, ...,\,m$, and $h^{(m+1)}_0\sspeq 1$. If $[d,a]\sspeq [1,0]$ the symbol $a_0 \sspeq 0$ can be omitted, and only the  $m$ symbols $a_j\sspeq a\sspp d\,j$ for $j\sspeq 1,\,2\, ...,\, m$ are active. For symmetric functions see \eg \cite {Krishnamurthy}, p. 53, and p. 54, eq $(46)$.
\psn 
The transition matrix property of the $\bf S2$[1,0] \sspeq  $\bf S2$ (see \cite{GKP}  p. 262, eq.\,$(6.10)$) generalizes to
\Beq
x^n\sspeq \sum_{m=0}^n\, \widehat{S2}(d,a;n,m)\,fallfac(d,a;x,m)\, ,\ \ n\sspin {\mathbb N}_0\, ,
\Eeq
where the generalized falling factorial is (see also Bala \cite{Bala} where this falling factorial appears in \Eq{15} as special $[t;d,0,c] _n$ in the signed $Stirling1$ context. See the present part C for the unsigned case) 
\Beq
fallfac(d,a;x,m)\sspdef \prod_{j=0}^{m-1} (x\sspm (a+j\,d)\, \ \ {\rm with}\ \ fallfac(d,a;x,0)\sspdef 1\, .
\Eeq
This can also be written in terms of the usual falling factorials $x^{\underline{n}}\sspdef \prod_{j=0}^{n-1}\,(x-j)$, for $n\sspin \mathbb N$ and $x^{\underline{0}}\sspdef 1$ as
\dstyle{fallfac(d,a;x,m)\sspeq d^m\,\left(\frac{x-a}{d}\right)^{\underline{m}}}.\psn
Using the binomial theorem in eq.\,$(11)$ and interchanging the sums shows that ${\bf S2}[d,a]$ (when the matrix elements are not specified we use this notation) can be written in terms of the usual $Stirling2$ numbers ${\bf S2}\sspeq {\bf S2}[1,0]$ as,
\Beq
S2(d,a;n,m)\sspeq \sum_{k=0}^n\,{\binomial{n}{k}}\,a^{n-k}\,d^k\,S2(k, m)\ .
\Eeq
For the inverse of this relation see {\sl Lemma 10}, \Eq{160}, in the proof section, part $C$.
\psn
A standard recurrence for {\sl Sheffer} row polynomials (\cite{Roman}, p. 50, Corollary 3.7.2) leads, with $PS2(d,a;n,x)\sspdef \sum_{m=0}^n\, S2(d,a;n,m)\,x^m$, to
\Beq
PS2(d,a;n,x)\sspeq [a\sspp d\,x\sspp d\,{\bf E}_x]\,PS2(d,a;n-1,x)\,, \ \ {\rm for}\ n \sspin \mathbb N\,, 
\Eeq 
with input $PS2(d,a;0,x)\sspeq 1$.\psn
The eq. $(8)$ version of the \ogf is not convenient to find $EPS(d,a;n,t)$ by inverse {\sl Laplace} transform because of the power $k\sspp 2$ instead of $k\sspp 1$. The solution is to consider first the \ogf of the powers (instead of the one of the sums  of powers) which can be found analogously to the $GPS$ case. Of course, if \Eq{8} has been proved the \ogf for the first difference sequence follows immediately. This will later lead to another form of $GPS$ which is amenable to find $EPS$.     
\Beqarray
GP(d,a;n,x)&\sspdef& \sum_{m=0}^{\infty}\, (a\sspp d\, m)^n\, x^m\,,\\
  &\sspeq & \sum_{k=0}^n\, S2(d,a;n,k)\, k!\,\frac{x^k}{(1-x)^{k+1}}\ .
\Eeqarray
From this the \egf can be computed directly based on \dstyle{{\cal L}^{-1}\left[ \frac{1}{(p-1)^{k+1}} \right]\sspeq  \frac{t^k}{k!}\, e^t}, using the linearity of the inverse {\sl Laplace} transform.
\Beq  
EP(d,a;n,t)\sspeq e^t\,\sum_{k=0}^n\,S2(d,a;n,k)\,t^k\,.
\Eeq
Continuing with the search for a more tractable form of $GPS(d,a;n,x)$  we apply another reordering identity on $GP(d,a;n,x)$, \viz
\Beq
\sum_{j=0}^n b^{(n)}_j\,\frac{x^j}{(1-x)^{j+1}}\sspeq \frac{1}{(1\sspm x)^{n+1}}\,\sum_{i=0}^n a^{(n)}_i\,x^i   \,,\ \ n\sspin \mathbb N_0\ ,
\Eeq
with
\Beqarray
a^{(n)}_i &\sspeq& \sum_{j=0}^i \, (-1)^{i-j}\, {\binomial{n-j}{i-j}}\, b^{(n)}_j,\\
b^{(n)}_j &\sspeq& \sum_{i=0}^{j} \, {\binomial{n-i}{j-i}}\, a^{(n)}_i\, .
\Eeqarray
Note that this reordering identity can not be applied to $GPS(d,a;n,x)$ because of the wrong power in the denominator. But here it produces
\Beqarray
GP(d,a;n,x) &\sspeq& \frac{1}{(1\sspm x)^{n+1}}\, PrEu(d,a;n,x)\,, \ {\rm with\ the\ polynomials} \\
PrEu(d,a;n,x) &\sspeq& \sum_{k=0}^n rEu(d,a;n,k)\,x^k\,, \ {\rm where} \\ 
rEu(d,a;n,k)&\sspeq& \sum_{j=0}^k (-1)^{k-j}\,{\binomial{n-j}{k-j}}\, S2(d,a;n,j)\,j!\, .  
\Eeqarray
Here $rEu(d,a;n,k)$ are generalized {\sl Euler}ian numbers, which constitute a number triangle (sometimes called {\sl Euler} triangle)  but compared with the usual {\sl Euler}ian triangle for $[d,a]\sspeq [1,0]$, given in {\sl Graham et al.} \cite{GKP}, Table 268, p. 268, or \seqnum{A173018}, the rows are reversed. The row reversed number triangle is shown in \seqnum{A123125}. This explains the $r$ in front of $Eu$ for {\sl Euler}ian. \pn
The inverse of the relation between ${\bf rEu}[d,a]$ and ${\bf S2}[d,a]$ is 
\Beq
S2(d.a;n,m)\,m! \sspeq \sum_{k=0}^m\,{\binomial{n-k}{m-k}}\, rEu(d,a;n,k)\, . 
\Eeq 
From the explicit form of ${\bf S2}[d,a]$ in eq.$(11)$ the one for ${\bf rEu}[d,a]$ follows by eq. $(23)$.
\Beq
rEu(d,a;n,k)\sspeq \sum_{j=0}^k\, (-1)^{k-j}\, {\binomial{n+1}{k-j}}\, (a+d\,j)^n\,.
\Eeq
In terms of the usual {\sl Eulerian} numbers one finds from eqs. $(28)$ with $(18)$ and eq. $(29)$ for $[d,a]\sspeq [1,0]$
\Beq
rEu(d,a;n,k) \sspeq \sum_{m=0}^n\,{\binomial{n}{m}}\,a^{n-m}\,d^m\,\sum_{p=0}^k\,(-1)^{k-p}\,{\binomial{n-m}{k-p}}\,rEu(m,p)\, .
\Eeq
Note that no formula analogous to \Eq{18} holds due to the different binomial structure in $rEu[d,a]$ of \Eq{28}.\psn
The three term recurrence for ${\bf rEu}[d,a]$ is
\Beq
rEu(d,a;n,m)\sspeq (d\,(n-m)\sspp (d-a))\, rEu(d,a;n-1,m-1) \,\sspp\, (a\sspp d\,m)\, rEu(d,a;n-1,m),   
\Eeq
for $n\sspgeq 1\,,\ m\sspeq 0,\,1,\,...,\,n$, with $rEu(d,a;n,-1)\sspeq 0$, $rEu(d,a;n,m)\sspeq 0$ for $n\sspkl m$ and $rEu(d,a;0,0)\sspeq 1$.\psn
The corresponding (ordinary, not exponential) row polynomials are
\Beq
PrEu(d,a;n,x)\sspdef \sum_{m=0}^n\,rEu(d,a;n,m)\,x^m\, ,\ \ n\sspin {\mathbb N}_0\,.
\Eeq
From eq.\,$(21)$ and eqs.\,$(26)$ with $(27)$ follows a relation between these ${\bf rEu}[d,a]$ row polynomials and those of the number triangle with entries $S2fac(d,a;n,m)\sspdef S2(d,a;n,m)\,m!$, named $PS2fac(d,a;n,x)$, \viz
\Beq
PrEu(d,a;n,x)\sspeq (1\sspm x)^n\,PS2fac\left(d,a;n,\frac{x}{1-x}\right)\, .
\Eeq
It may be noted, in passing, that the transformation \dstyle{y\sspeq \frac{x}{1\sspp x}}, or \dstyle{x\sspeq \frac{y}{1\sspm x}} is called {\sl Euler}'s transformation (see, \eg \cite{Hardy}, p. 191, last row).
\psn
From this preceding relation the \egf (\ie the e.g.f. for the row reversed {\sl Euler}ian triangle) follows:
\Beq
{\fbox{\color{cyan}$EPrEu(d,a;t,x)$}}\sspeq \frac{(1\sspm x) \,e^{a\,(1-x)\,t}}{1\sspm x\,e^{d\,(1-x)\,t}}\, .
\Eeq
This is not a {\sl Sheffer} structure, not even one of the more general {\sl Brenke} type $g(z)\, B(x\,z)$, \cite{Brenke}, \cite{Chihara}, p. 167.\psn
The \egf of the row sums ($x\sspto 1$) of ${\bf rEu}[d,a]$ is obtained \via\ {\sl l'H\^opital}'s rule as \dstyle{\frac{1}{1\sspm d\,t}}, independently of $a$.\psn 
A one parameter $k$-family of generalized Eulerian polynomials $A_{n,k}(x)$ with coefficient triangles has been considered by {\sl Luschny} \cite{Luschny2}. The coefficients of $A_{n,k}(x)$ build ${\bf Eu}[k,1] \sspeq {\bf rEu}[k,k-1]$.\psn Now the new form of $GPS(d,a;n,x)$ is simply obtained by multiplying $GPS(d,a;n,x)$ with \dstyle{\frac{1}{1\sspm x}} because this is the rule to obtain the \ogf for partial sums of a sequence from the \ogf of the sequence.
\Beq
GPS(d,a;n,x) \sspeq \frac{1}{(1\sspm x)^{n+2}}\, PrEu(d,a;n,x)\, . 
\Eeq
There is still this power $n+2$ but now we can use the reordering identity, eq.\,$(23)$ with $n$ replaced by $n+1$: 
\Beq
\sum_{j=0}^{n+1}\, b^{(n+1)}_j\, {\frac{x^j}{(1\sspm x)^{j+1}}}\sspeq \frac{1}{(1\sspm x)^{(n+1)+1}} \sum_{i=0}^{n+1}\,a_i^{(n+1)}\,x^i\, , 
\Eeq
with (see eq.\,$(30)$) 
\Beq
a^{(n+1)}_i\sspeq rEu(d,a;n,i)\sspeq \sum_{p=0}^i\,(-1)^{i-p}\,{\binomial{n+1}{i-p}}\, (a\sspp d\,p)^n\, . 
\Eeq
Note that $a^{(n+1)}_{n+1}\sspeq 0$ because $PrEu(d,a;n,x)$ has degree $n$. \psn
Note that now the $x$ dependence is amenable for a later inverse {\sl Laplace} transform. The calculation of $b^{(n+1)}_j$ is a bit lengthy but it turns out to have a nice form (we add the $[d,a]$ parameters).
\Beq
b_j^{(n+1)}(d,a)\sspeq S2(d,a;n,j)\,j! \sspp S2(d,a;n,j-1)\,(j-1)! \sspfed \Sigma S2(d,a;n,j)\, ,
\Eeq
leading finally to the result for the \egf
\Beq
{\fbox{\color{magenta}$EPS(d,a;n,t)$}} \sspeq e^t\, \Sigma S2(d,a;n,j) \frac{t^j}{j!}\ , n\sspin {\mathbb N}_0\, .
\Eeq
Let us recapitulate the detour we made in a diagram referring to eq.\,$(23)$ for obtaining two versions of $GP$ or $GPS$: 
\Beqarray
{\bf GPSv1}\ &{\buildrel {(23)} \over \nrightarrow}&\  {\bf GPSv2}\,\ \ {\buildrel {\cal L}^{-1} \over \rightarrow}\ \ {\bf EPS}\,  \nonumber \\
\downarrow \phantom{xx} &&\hskip .5cm \uparrow \nonumber\\
{\bf GPv1}\  &{\buildrel {(23)} \over \rightarrow}&\ \, {\bf GPv2}\, .
\Eeqarray
\pbn
{\bf B) Generalized Faulhaber formula and Bernoulli polynomials}\pn
The next topic is to find for the power sum $PS(d,a;n,m)$ a formula in terms of {\sl Bernoulli} polynomials evaluated appropriately. This formula has bee named {\sl Faulhaber} formula for the ordinary $[d,a]\sspeq [1,0]$ case by {\sl Conway} and {\sl Guy} \cite{CoGuy}, p. 106. {\sl Faulhaber} used the numbers, later called {\sl Bernoulli} numbers  by {\sl de Moivre} and {\sl Euler} (see \cite{Edwards1}, \cite{Edwards2}), already by 1631 before {\sl Jakob I Bernoulli}. For this formula see \cite{GKP}, p. 367. eq. (7.79), and \cite{Koecher}, p. 167 eq. (1). Here it is, with our definition of $PS(n,m) = PS(1,0;n,m)$,
\Beq
PS(n,m) \sspeq \delta_{n,0} \sspp \frac{1}{n+1}\, \left( B(n+1,\, x=m+1)\sspm B(n+1,\,x=1)\right)\, , \ \ n\sspin \mathbb N_0\, , 
\Eeq  
where $\delta_{n,0}\sspeq [n\sspeq 0]$ is the {\sl Kronecker} symbol: $1$ if $n\sspeq 0$ and 0 otherwise. The {\sl Bernoulli} numbers are defined recursively by (see \cite{GKP}, p. 284, eq. $(6.79)$)
\Beq
B(n)\sspdef \frac{1}{n\sspp 1}\,\left(\delta_{n,0}\sspm \sum_{k=0}^{n-1}\, {\binomial{n+1}{k}} B(k) \right)\, \ {\rm for}\ \ n\sspin \mathbb N, \ {\rm with}\ B(0)\sspeq 1\, .  
\Eeq
They have \dstyle{B(1)\sspeq -\frac{1}{2}} and are found in OEIS \cite{OEIS} under \seqnum{A027641}\, /\,\seqnum{A027642}.
The corresponding {\sl Bernoulli} polynomials are 
\Beq
B(n,\,x)\sspdef \sum_{m=0}^m {\binomial{n}{m}} B(n-m)\,x^m\,. 
\Eeq  
Their coefficient tables are given in \seqnum{A196838}\,/\,\seqnum{A196839} or \seqnum{A053382}\, /\, \seqnum{A053383} for rising or falling powers of $x$, respectively.\psn
For the generalized case one finds for the power sums $PS(d,a;n,m)$ from eqs. $(5)$ and $(42)$ the following {\sl Faulhaber} formula in terms of ordinary {\sl Bernoulli} polynomials 
\Beq
{\fbox{\color{blue}$PS(d,a;n,m)$}}\sspeq \sum_{k=0}^n\,{\binomial{n}{k}}\, a^{n-k}\, d^k\, \left [ \delta_{k,0} \sspp \frac{1}{k+1}\,\left(B(k+1,\,x=m+1)\sspm B(k+1,\,x=1) \right) \right]\,.
\Eeq
\pbn
But the idea is to find the analogon of formula $(42)$ with generalized {\sl Bernoulli} polynomials.\psn
An obvious generalization of the {\sl Bernoulli} numbers is 
\Beq
B(d,a;n) \sspdef \sum_{m=0}^n\,(-1)^m\,\frac{1}{m+1}\,S2(d,a;n,m)\,m!\,,\ n\sspin {\mathbb N}_0\,.
\Eeq
For $[d,a]\sspeq [1,0]$ see, \eg {\sl Charalambides} \cite{Charalambides}, or the formula and Maple section of \seqnum{A027641}.
\psn
From eq.$(18)$ one finds $B[d,a]$ in terms of $B$.
\Beq
B(d,a;n) \sspdef \sum_{m=0}^n {\binomial{n}{m}}\,a^{n-m}\,d^m\,B(m)\, .
\Eeq
The \egf of $\{B(d,a;n)\}_{n=0}^{\infty}$ is
\Beq
EB(d,a;t)\sspeq \frac{d\,t\,e^{a\,t}}{e^{d\,t}\sspm 1}\, .
\Eeq
The corresponding generalized {\sl Bernoulli} polynomials are (compare with eq.\,$(44)$)
\Beq
B(d,a;n,x) \sspdef \sum_{m=0}^n {\binomial{n}{m}}\,B(d,a;n-m)\,x^m \,,
\Eeq
and from eq.\,$(47)$ they can also be written in terms of $\{B(m)\}_{m=0}^n$ as
\Beq
B(d,a;n,x) \sspdef \sum_{m=0}^n {\binomial{n}{m}}\,d^m\,B(m)\,(a\sspp x)^{n-m} \,.
\Eeq
Their \egf is, either from \Eq{49} or $(50)$,
\Beq
EB(d,a;t,x)\sspeq \frac{d\,t\,e^{a\,t}}{e^{d\,t}\sspm 1}\,e^{x\,t}\, ,
\Eeq
identifying their coefficients as {\sl Sheffer} arrays \dstyle{\left(\frac{d\,z\,e^{a\,z}}  {e^{d\,z}\sspm 1},\,z \right)}. Such arrays are of the so called {\sl Appell} type (compare {\sl Roman} \cite{Roman}, pp. 26 - 28, with a different notation).
It turns out that in order to obtain a generalized {\sl Faulhaber} formula in terms of {\sl Bernoulli} polynomials the $B[d,a]$ just introduced are not quite the ones needed. In fact, they are too general. One has to work with the polynomials depending only on $d$, \viz
\Beq
{\fbox{\color{bittersweet}$B(d;n,x)$}}\sspeq \sum_{m=0}^n {\binomial{n}{m}}\, B(d;n-m)\, x^m\, ,
\Eeq
with the generalized {\sl Bernoulli} numbers
\Beq
B(d;n)\sspdef B(d,a=0;n) = d^n\,B(n)\,, \  n\sspin {\mathbb N}_0\, .
\Eeq
They can also be obtained by exponential convolution of the more general ones with the sequence $\{-a^n\}_{n=1}^{\infty}$.
\Beq
B(d;n)\sspeq \sum_{m=0}^n {\binomial{n}{m}}\, B(d,a;n-m)\, (-a)^m\,.
\Eeq
For $a\sspeq 0$ only $m\sspeq 0$ survives and B(d,0;n) results. But also for non-vanishing $a$ the $a$ dependence drops out, as can be seen from the \egf of the sequence on the {\it r.h.s.}, using eq. $(48)$.
\Beq
\frac{d\,t\,e^{a\,t}}{e^{d\,t}\sspm 1}\, e^{-a\,t} \sspeq  \frac{d\,t}{e^{d\,t}\sspm 1}\sspeq EB(d,a=0;t)\sspfed EB(d;t)\,.
\Eeq
For $B(d;n)$ with $d\sspeq 2,\,3$ and $4$ see $(-1)^n$\seqnum{A239275}$(n)$/\seqnum{A141459}$(n)$,\ \seqnum{A285863}$(n)$/\seqnum{A285068}$(n)$ and\pn
\seqnum{A288873}$(n)$/\seqnum{A141459}$(n)$.\psn
The \egf of the polynomial system $\{B(d;n,x)\}_{n=0}^{\infty}$ of eq.\,$(52)$ is
\Beq
 EB(d;t,x)\sspeq  \frac{d\,t}{e^{d\,t}\sspm 1}\, e^{x\,t}\sspeq EB(d,a=0;t,x)\,.
\Eeq
The {\sl Appell} type {\sl Sheffer} structure is obvious.\psn
Now the stage is set for giving the result for the generalized {\sl Faulhaber} formula in terms of the polynomials $B(d,n,x)$.
\Beqarray
{\fbox{{\color{blue}$PS(d,a;n,m)$}}}&\sspeq& \frac{1}{d\,(n+1)}\,\left[B(d;n+1,x \sspeq a\sspp d\,(m+1))\sspm B(d;n+1,x \sspeq d) \right. \nonumber \\
&&\hskip 1.5cm \sspm \left. B(d;n+1,x \sspeq a)\sspp B(d,n+1,x=0) \sspp d\,\delta_{n,0}\right]\,. 
\Eeqarray 
Here $B(d,n+1,x=0) = B(d;n+1)\sspeq d^{n+1}\,B(n+1)$, and the {\sl Kronecker} symbol enters because of our definition of $PS(d,a;n,m)$ where the sum starts with $j=0$, not with $1$.
\psn
The generalized {\sl Lah} numbers ${\bf L}[d,a]$ are discussed in the proof section 2, C) 4.\pbn\psn
{\bf C) Generalized Stirling1 numbers}\psn
As elements of the {\sl Sheffer} group the inverse of the (infinite, lower triangular) matrix $\bf S2[d,a]$ exists and is called $\bf S1[d,a]$. This is therefore a generalized {\sl Stirling} number triangle of the first kind. 
\Beq
 {\bf S2}[d,a]\cdot {\bf S1}[d,a]\sspeq {\bf 1} \sspeq {\bf S1}[d,a]\cdot {\bf S2}[d,a]\,,
\Eeq
with the (infinite dimensional) identity matrix {\bf 1}. For practical purposes it is sufficient to consider the  finite dimensional case of $N\ssptimes N$ matrices. ${\bf S1}[d,a]$ is a signed matrix with fractional entries for $d\sspneq 1$. Therefore, in order to have non-negative entries one considers ${\bf S1p}[d,a]$ with entries $S1p(d,a;n,m) \sspdef (-1)^{n-m}\,S1(d,a;n,m)$. But in the combinatorial context also a scaling is needed to obtain a non-negative integer matrix ${\bf \widehat{S1p}}[d,a]$ with diagonal entries $1$ (\ie monic row polynomials). This is done by scaling the ${\bf S1p}[d,a]$ rows $n$ with $d^n$. \psn
We then have the {\sl Sheffer} structures
\Beq
{\bf S1}[d,a]\sspeq \left(\frac{1}{(1\sspp x)^{\frac{a}{d}}},\, \frac{1}{d}\,log(1\sspp x)\right)\,,\hskip 1cm \rm {and} \hskip 1cm 
{\fbox{\color{bleudefrance}${\bf\widehat{S1p}}[d,a]$}}\sspeq \left(\frac{1}{(1\sspm d\,x)^{\frac{a}{d}}},\, -\frac{1}{d}\,log(1\sspm d\,x)\right)\,. 
\Eeq
The ${\bf \widehat {S2}}[d,a]$ matrices (see \Eq{15}) which have scaled matrix elements $\widehat{S2}(d,a;n,m)\sspeq S2(d,a;n,m)/d^m$ have the signed inverse matrices \dstyle{{\bf \widehat{S1}}[d,a]\sspeq \left((1\sspp d\,x)^{-\frac{a}{d}},\,\frac{1}{d}\,\log(1\sspp d\,x)\right)}.\psn
The signed ${\bf \widehat {S1}}[d,a]$ matrices have been considered by {\sl Bala} \cite{Bala} as $s_{(d,0,a)}$. In {\sl Luschny} \cite{Luschny1} the $SF-C$ matrices are our ${\bf \widehat{S1p}}[d,d-1]$, and the $SF-CS$ matrices are the unsigned inverse matrices of  ${\bf \widehat {S2}}[d,d-1]$.\psn 
The {\sl Sheffer} structure of ${\bf \widehat{S1p}}[d,a]$ means that the \egf of column $m$ is
\Beq
{\fbox{\color{bleudefrance}$E\widehat{S1p}Col(d,a;t,m)$}}\sspeq \frac{1}{(1\sspm d\,t)^{\frac{a}{d}}}\, \frac{1}{m!}\,\left(-\frac{1}{d}\,log(1\sspm d\,t)\right)^m \,, \ \ m\sspin {\mathbb N}_0\ .
\Eeq
There seems not to exist a simple form for the corresponding \ogf.\psn
The three term recurrence for the ${\bf \widehat{S1p}}[d,a]$ matrix entries is 
\Beq
\widehat{S1p}(d,a;n,m) \sspeq \widehat{S1p}(d,a;n-1,m-1) \sspp (d\,n\sspm (d-a))\, \widehat{S1p}(d,a;n-1,m), \ \ {\rm for}\ n\sspgeq 1\,,\ m\sspeq 0,\,1,\,...,\,n\,,  
\Eeq
with $\widehat {S1p}(d,a;n,-1)\sspeq 0$, $\widehat {S1p}(d,a;n,m)\sspeq 0$ for $n\sspkl m$ and $\widehat {S1p}(d,a;0,0)\sspeq 1$.\psn
The usual transition from the monomial basis $\{x^n\}_{n=0}^{\infty}$ to the rising factorials (see \cite{GKP}, p. 263, eq.\,$(6.11)$) generalizes to the following identification of the row polynomials of ${\bf \widehat{S1p}}[da]$
\Beq
 {\fbox{\color{bleudefrance}$P\widehat{S1p}(d,a;n,x)$}}\sspdef \sum_{m=0}^n\,\widehat{S1p}(d,a;n,m)\,x^m\sspeq risefac(d,a;x,n)\, ,
\Eeq
with the generalized rising factorials (compare this with the generalized falling factorials eq.\,$(17) $)
\Beq
risefac(d,a;x,n)\sspdef \prod_{j=0}^{n-1} (x\sspp (a+j\,d))\, \ \ {\rm with}\ \ risefac(d,a;x,0)\sspdef 1 .
\Eeq
This can be rewritten also in terms of the usual rising factorial $x^{\overline{n}}\sspdef \prod_{j=0}^n\,(x+j)$ for $n\sspin \mathbb N$ and $x^{\overline{0}}\sspdef 1$ as \dstyle{risefac(d,a;x,n)\sspeq d^n\,\left(\frac{x+a}{d} \right)^{\overline{n}}}. In terms of falling factorials this is  \dstyle{risefac(d,a;x,n)\sspeq (-d)^n\,\left(\frac{-(x+a)}{d} \right)^{\underline{n}}}.\psn
This identification implies \via\ {\sl Vieta}'s theorem that the coefficients of the monic polynomial $P\widehat{S1p}(d,a;n,x)$ are the elementary symmetric functions $\sigma^{(n)}_{n-m}(a_0,a_1,...,a_{n-1})$ in the indeterminates $\{a_j\}_{j=0}^{n-1}$ given by $a_j\sspdef a\sspp j\,d$, with $\sigma^{n}_{0}\sspdef 1$. Sometimes \dstyle{\sigma^{(n)}_{n-m}[d,a]} is used for these symmetric functions.  Thus
\Beq
\widehat{S1p}(d,a;n,m)\sspeq \sigma^{(n)}_{n-m}(a_0,a_1,...,a_{n-1}),\ {\rm with}\  a_j\sspeq a\sspp j\,d\,. 
\Eeq 
If $d\sspeq 1$ (and $a\sspeq 0$)  $a_0\sspeq 0$ does not contribute and one can write $\widehat{S1p}(1,0;n,m)\sspeq S1p(n,\, m)\sspeq \sigma^{(n-1)}_{n-m}(1,\,2,\,...,\, n-1)$.\psn 
Sorting in falling powers of $a$ one obtains the formula for $\widehat{S1p}(d,a;n,m)$ in terms of the usual unsigned {\sl Stirling1} numbers $S1p(n,\,m) =$\seqnum{A132393}$(n,\,m)$.
\Beq
\widehat{S1p}(d,a;n,m)\sspeq \sum_{j=0}^{n-m}\,{\binomial{n-j}{m}}\,S1p(n,\,n-j)\, a^{n-m-j}\, d^j\  \sspeq \sum_{j=m}^{n}\,{\binomial{j}{m}}\,S1p(n,\,j)\, a^{j-m}\, d^{n-j} \, . 
\Eeq
This satisfies the recurrence relation \Eq{61}.\psn
The standard {\sl Sheffer} recurrence (\cite{Roman}, p. 50, Corollary 3.7.2) for these row polynomials boils down to \Beq
P\widehat{S1p}(d,a;n,x)\sspeq (x+a)\,P\widehat{S1p}(d,a;n\sspm 1,x\sspp d)\,, \ \  n \sspin \mathbb N\,,
\Eeq
with input $P\widehat{S1p}(d,a;0,x)\sspeq 1 $.\psn
From the {\sl Sheffer} property the \egf of these row polynomials, \ie the \egf of the number triangle ${\bf \widehat{S1p}}[d,a]$, is 
\Beq
{\fbox{\color{bleudefrance}$EP\widehat{S1p}(d,a;t,x)$}}\sspeq \frac{1}{(1\sspm d\,t)^{\frac{a}{d}}}\,exp\left(-x\,\frac{1}{d}\,log(1\sspm d\,t)\right)\sspeq \frac{1}{(1\sspm d\,t)^{\frac{a+x}{d}}}\ .
\Eeq
For the {\sl Meixer} type recurrence see the proof section 2, C), 7.\psn
A more involved problem is to find the generalization of the formula giving $\widehat{S1p}(d,a;n,m)$ in terms of the column scaled $\widehat{S2}(d,a;n,m)$ elements. The standard {\sl Schl\"omilch} formula is (see, e.g., \cite{Charalambides}, p. 290, eq. $(8.20)$ for the signed $S1$ entries) 
\Beq
S1p(n,\,m)\sspeq (-1)^{n-m}\,\sum_{k=0}^{n-m}\, (-1)^k\,{\binomial{n+k-1}{m-1}}\, {\binomial{2\,n - m}{n-m-k}}\, S2(n-m+k,\,k)\ .
\Eeq
The direct proof starts with \Eq{65}. Inserting the {\sl Schl\"omilch} formula just given, then using the inverse of \Eq{18} leads to
\Beqarray
{\fbox{\color{bleudefrance}$\widehat{S1p}(d,a;n,m)$}}&\sspeq& a^{n-m}\,\sum_{j=m}^n\,{\binomial{j}{m}}\,\sum_{k=0}^{n-j}\,{\binomial{n-k-1}{j-1}}\,{\binomial{2\,n-j}{n-j-k}}\,a^k\, * \nonumber \\
&& *\sum_{l=0}^{n-j+k}\,(-1)^l\,{\binomial{n-j+k}{l}}\,a^{-l}\,\widehat{S2}(d,a;l,k)\,,\ {\rm for}\  n\sspgeq m\sspgeq 0\,.   
\Eeqarray
Note that this result also holds for $a\sspeq 0$ because then $a^{n-m+k-l}$ becomes $\delta_{0,n-m+k-l}$ (from $0^0\sspeq 1$) leading to a collapse of the $l-$sum, and the remaining two sums produce $\widehat{S1p}(d,0;n,m)\sspeq d^{n-m}\,S1p(n,\,m)$. \psn
Also the known result for $m\sspeq 0$ from \Eq{62} is recovered, \viz\ $\widehat{S1p}(d,a;n,0)\sspeq risefac(d,a;0,n)\sspeq d^n\,\left(\frac{a}{d}\right)^{\overline{n}}$.\psn
For another formula, following from a proof along the lines of the ordinary formula in \cite{Charalambides}, p. 290, see the proof section {\it C}, \Eq{175}. 
\psn
The inverse of the generalized {\sl Lah} matrix ${\bf L}^{-1}[d,a]$ is discussed in the proof section 2, C) 4.\pbn\psn
{\bf D) Combinatorial Interpretation }\psn
{\bf I)} ${\bf \widehat{S2}}[d,a]$\psn
The \ogf, eq. $(14)$, divided by $d^m$, which generates the complete homogeneous symmetric functions $h^{(m+1)}_{n-m}[d,a]$ of degree $n-m$ of the $m+1$ symbols $a_j\sspeq a\sspp d\,j$, $j\sspeq 0,\,1,\,...,\, m$, leads immediately to the following combinatorial interpretation of \dstyle{\widehat{S2}(d,a;n,m)\sspdef \frac{1}{d^m}\, S2(d,a;n,m)} (see \Eq{15}).\psn
$\widehat{S2}(d,a;n,m)$ is for $d\sspgeq 2$ the (dimensionless) total volume of the $multichoose(m+1,\, n-m)\sspeq {\binomial{n}{m}}$  hyper-cubes and hyper-cuboids (polytopes) of dimension $n-m$ which are build from the $n-m$ orthogonal $\mathbb Z^{n-m}$ vectors of lengths taken from the repertoire $a_j\sspeq a\sspp d\,j$, $j\sspeq  0,\,1,\,...,\, m$. \psn
For $d\sspeq 1$ (and $a\sspeq 0$), the standard {\sl Stirling2} case ${\bf S2}\sspeq {\bf S2}[1,0]\sspeq {\widehat {\bf S2}}[1,0]$, $a_0\sspeq 0$ does not contribute and the $n-m$ vectors are from the set $\{1,\,2,\, ...,\, m\}$ for the $multichoose(m,\, n-m)\sspeq {\binomial{n-1}{m-1}}$ polytopes.\psn
Some examples:\psn
a) ${\widehat S2}(1,0;3,2)\sspeq S2(3,2)\sspeq 3$ from the \dstyle{{\binomial{2}{1}}\sspeq 2} polytopes of dimension $1$ with basis lengths  $1,\,2$, \ie two lines of length $1$ and $2$ with total length $3$.\psn
b) ${\widehat S2}(2,1;3,2)\sspeq 9$ (see \seqnum{A039755})  from the \dstyle{{\binomial{3}{2}}\sspeq 3} polytopes of dimension $1$ with basis lengths  $1,\,3,\, 5$, \ie three lines of total length $9$.\psn
c)  ${\widehat S2}(3,2;3,1)\sspeq 39$ (see \seqnum{A225468}) from the \dstyle{{\binomial{3}{1}}\sspeq 3} polytopes of dimension $2$ with basis lengths from the set $\{2,\,5\}$, \ie  two squares of area $2^2$ and $5^2$ and a rectangle of area $2^1\,5^1$,  giving total area $4\sspp 25\sspp 10\sspeq 39$. \pbn
{\bf II)} ${\bf \widehat{S1p}}[d,a]$\psn
From eq. ($62$) for the row polynomials and the implied elementary symmetric function formula for $\widehat{S1p}(d,a;n,m)$ of eq.\,$(64)$ one has the combinatorial interpretation.\psn
$\widehat{S1p}(d,a;n,m)$ is  for $d\sspgeq 2$ the total volume of the ${\binomial{n}{n-m}}$ hyper-cuboids of dimension $n-m$ formed from the $n-m$ orthogonal $\mathbb Z^{n-m}$ vectors with distinct (dimensionless) lengths from the $n$-set $\{a\sspp d\,j\,|\, j=0,\,1,\,...,\, n-1\}$. For $[d,a]\sspeq [1,0]$ the ordinary unsigned {\sl Stirling1} number $S1p(n,\, m)$ (see \seqnum{A132393}) gives the total volume of the ${\binomial{n-1}{n-m}}$ hyper-cuboids of dimension $n-m$ formed from the $n-m$ orthogonal $\mathbb Z^{n-m}$ vectors with distinct (dimensionless) lengths from the $(n-1)$-set $\{1,\,2,\,...\, n-1\}$.\psn
Some examples:\psn
a) $\widehat{S1p}(1,0;4,2)\sspeq S1p(4,\, 2)\sspeq 11$ (see \seqnum{A132393}) from the \dstyle{{\binomial{3}{2}}\sspeq 3} hyper-cuboids of dimension $2$ with distinct basis vector lengths from the set $\{1,\,2,\, 3\}$, \ie six rectangles of area $1\sspcdot 2$, $1\sspcdot 3$ and $2 \sspcdot 3$, with total area $2\sspp 3\sspp 6\sspeq 11$.\psn
b) $\widehat{S1p}(2,1;4,1)\sspeq 176$ (see \seqnum{A028338}) from the  \dstyle{{\binomial{4}{3}}\sspeq 4} hyper-cuboids of dimension $3$ with distinct basis vector lengths from the set $\{1,\,3,\,5,\,7\}$, \ie four cuboids with volumes $1\sspcdot 3 \sspcdot 5$, $1\sspcdot 3\sspcdot 7$, $1 \sspcdot 5\sspcdot 7$ and $3\sspcdot 5 \sspcdot 7$, adding to $176$.\psn
c) $\widehat{S1p}(3,1;4,0)\sspeq 280$ (see \seqnum{A286718}) from the  \dstyle{{\binomial{4}{4}}\sspeq 1} hyper-cuboid of dimension $4$
with distinct basis vector lengths from the set $\{1,\,4,\,7,\,10\}$, \ie the $4D$ hyper-cuboid with volume $1\sspcdot 4 \sspcdot 7 \sspcdot 10\sspeq 280$. 
\pbn
Two remarks: The first column sequences $\{\widehat{S1p}(d,1;n,0)\}_{n=1}^{\infty}$ have also an interpretation as numbers of $(d+1)$-ary rooted increasing trees with $n$ vertices, including the root vertex. This is the sequence $\{S(k=d+1;n,1)\}$ of generalized {\sl Stirling2} numbers with parameter $k$ in the notation of \cite{WLang}, eq. ($5$). The reason is the \egf called there $g2(k=d+1;x)\sspeq -1 \sspp (1 \sspp (1-(d+1))\,x)^{\frac{1}{1-(d+1)}} \sspeq -1 \sspp (1\sspm d\,x)^{-\frac{1}{d}}$ which is the \egf of the $m=0$ column of ${\bf \widehat{S1p}}[d,1]$, viz $E\widehat{S1p}Col(d,1;x,0)\sspeq (1\sspm d)^{-\frac{1}{d}}$ (see \Eq{60}) but with the $n\sspeq 0$ entry removed. See the instances \seqnum{A001147}, \seqnum{A007559}, \seqnum{A007696}, \seqnum{A008548}, ... for $d\sspeq 2,\, 3,\, 4,\, 5\, ...$, respectively. They are, for $n >=1 $, the number of $3,\,4,\,5,\,6,\, ...$-ary increasing rooted trees. \psn
Similarly, the first column sequences $\{\widehat{S1p}(d,d-1;n,0)\}_{n=1}^{\infty}$ are related to another variety of increasing trees given by the \egf $g2p(k=d-1;x)\sspeq 1\sspm (1\sspm  d\,x)^{\frac{1}{d}}$ for the sequence $\{|S(-k=1-d;n,1)|\}_{n=0}^{\infty}$ of \cite{WLang}, eq.\,($6$). This is related to the \egf of the sequence $\{\widehat{S1p}(d,d-1,n,0)\}_{n=0}^{\infty}$, \ie $E\widehat{S1p}Col(d,d-1;x,0)\sspeq  (1\sspm d\,x)^{-\frac{d-1}{d}}$ by integrating and adding $1$: $\int dx\,g\widehat{S1p}(d,d-1;x) \sspp 1\sspeq g2p(d-1;x)$. \pn
Some instances are: \seqnum{A001147}, \seqnum{A008544}, \seqnum{A008545}, \seqnum{A008546}, ... for $d\sspeq 2,\, 3,\, 4,\, 5\, ...$, respectively. 
\pbn
The combinatorial interpretation of ${\bf rEu}[da]$ should also be consodered.
\pbn
\vskip 1cm
\section{Proofs} 
\hskip 1cm In the following proofs the setting of formal power series is used. No convergence issues are considered. Infinite sums are interchanged (one could use alternatively a large cutoff). Differentiation as well as integration will also be interchanged with infinite sums. Only statements which are not already obvious from the main text are proved here. Note that for binomials the definition of \cite{GKP}, p. 154, eq.\,$(5.1)$ is taken. This is not the definition used by {\sl Maple13} \cite{Maple}. Also $0^0\sspdef 1$. The symmetry of binomial coefficients is used {\sl ad libitum} (but the upper number in the binomial has to be a non-negative integer).\psn
{\bf A) Proofs of section 1\,A}\psn
{\bf 1. Proof of eqs. $\bf (8)$ to $\bf (11)$}\psn
{\bf Lemma 1:}  
With the notation of eq. ($9$):
\Beq
(a\,{\bf 1} \sspp d\,{\bf E}_x)^n\,x^j\sspeq (a\sspp d\,j)^n\, x^j\ .
\Eeq
{\bf Proof:} Trivial, by induction over $n\sspin \mathbb N_0$ with $j\sspin \mathbb N_0$.\psn
For the \ogf $GPS(d,a;n,x)$ from eq.\,($2$) with eq.\,($1$) one has, after an exchange of the two sums, inserting $x^j\,x^{-j}$ and application of {\sl Lemma 1}: 
\Beqarray
GPS(d,a;n,x) &\sspeq& \sum_{m=0}^{\infty}\,x^m\, \sum_{j=0}^m\, (a\sspp d\,j)^n \sspeq \sum_{j=0}^{\infty}\, (a\sspp d\,j)^n\,\sum_{m=j}^{\infty}\,x^m \nonumber \\
&\sspeq& \sum_{j=0}^{\infty}\, \left( (a\,{\bf 1} \sspp d\,{\bf E}_x)^n\,x^j\right)\, \sum_{m=j}^{\infty}\, x^{m-j}\sspeq  \sum_{j=0}^{\infty}\, \left( (a\,{\bf 1} \sspp d\,{\bf E}_x)^n\,x^j\right)\,\frac{1}{1\sspm x}\, . 
\Eeqarray
After summing over $j$, the reordering of differentials from eq.\,($9$), \ie the definition of $S2(d,a;n,m)$, is used:
\Beqarray
\hskip 1.3cm &\sspeq& (a\,{\bf 1} \sspp d\, {\bf E}_x)^n\,\frac{1}{(1\sspm x)^2}\sspeq \sum_{k=0}^n\,S2(d,a;n,k)\,x^k\,{\bf d}_x^k\,\frac{1}{(1\sspm x)^2} \nonumber \\
&\sspeq& \sum_{k=0}^n\,S2(d,a;n,k)\,k!\,\frac{x^k}{(1\sspm x)^{2+k}}\, .
\Eeqarray
The three term recurrence eq.\,($10$) of the number triangle  $\{S2(d,a;n,k)\}$ follows from the definition eq. ($9$):
\Beqarray
\sum_{m=0}^{n+1}\, S2(d,a;n+1,m)\,x^m\,{\bf d}_x^m &\sspeq& (a\,{\bf 1} \sspp d\, {\bf E}_x)\,(a\,{\bf 1} \sspp d\, {\bf E}_x)^n \sspeq (a\,{\bf 1} \sspp d\,{\bf E}_x)\,\sum_{m=0}^{n}\, S2(d,a;n,m)x^m\,{\bf d}_x^m \nonumber \\
&\sspeq & \sum_{m=0}^{n}\, S2(d,a;n,m)\,a\, x^m\,{\bf d}_x^m \sspp   \sum_{m=0}^{n}\, S2(d,a;n,m)\,d\,\left(m\, x^m\,{\bf d}_x^m \sspp x^{m+1}\,{\bf d}_x ^{m+1}\right) \nonumber \\ 
&\sspeq& \sum_{m=0}^{n+1}\, S2(d,a;n,m)\,(a\, \sspp d\,m)\,x^m\,{\bf d}_x^m \sspp \sum_{m=1}^{n+1}\, S2(d,a;n,m-1)\,d\,x^m\,{\bf d}_x^m\ . 
\Eeqarray
In the second to last sum the $m\sspeq n+1$ term has been added due to the triangle condition $S2(d,a;n,m)\sspeq 0$ if $n\sspkl m$. In the last sum the lower index $m\sspeq 0$ can be added because of the input condition $S2(d,a;n,-1)\sspeq 0$. Comparing powers of $x^m\,{\bf d}_x^m$ then leads to the recurrence eq.\,($10$) after the change $n\sspto n-1$.\pbn
The explicit form of $S(d,a;n,m)$ from \Eq{11} satisfies the recurrence \Eq{10} together with the inputs because
\Beqarray
&& d\, S2(d,a;n-1,m-1) \sspp (a\sspp d\,m)\,S2(d,a;n-1,m) \nonumber \\
&&\sspeq \sum_{k=0}^{m-1}\frac{d}{(m-1)!}\,(-1)^{m-k}\,(-1)\,{\binomial{m-1}{k}}\,(a\sspp d\,k)^{n-1} \sspp \sum_{k=0}^{m}\frac{1}{m!}\,(-1)^{m-k}\,{\binomial{m}{k}}\,(a\sspp d\,m)\,(a\sspp d\,k)^{n-1} \nonumber \\
&\sspeq& \frac{1}{m!}\,\sum_{k=0}^{m}\,(-1)^{m-k}\,{\binomial{m}{k}}\,\left[-d\,m \frac{{\binomial{m-1}{k}}}{{\binomial{m}{k}}} \sspp (a\sspp d\,m)\right]\,(a\sspp d\,k)^{n-1} \, .
\Eeqarray
I the first term of the last sum the term $k\sspeq m$ does not contribute because of the binomial. Now the term within the bracket becomes
\Beq
\left[ ...\right]\sspeq -d\,\frac{m}{m}\,(m\sspm k) \sspp (a\sspp d\,m)\sspeq a\sspp d\,k\,,
\Eeq
leading to the \lhs of the recurrence, viz\, $S2(d,a;n,m)$ of eq.\,($11$)
\psn
For the instances of $S2[d,a]$ for $[d,a]\speq [1,0],\,[2,1],\,[3,1],\, [3,2],\,[4,1],\, ]4,3]$ see \seqnum{A048993},\ \seqnum{A154537},\, \seqnum{A282629},\, \seqnum{A225466},\, \seqnum{A285061},\, \seqnum{A225467}, respectively. \pbn
\psn
We give also the recurrence for $S2fac(d,a;n,m)\spdef S2(d,a;n,m)\,m!$:
\Beq
S2fac(d,a;n,m)\sspeq m\,d\, S2fac(d,a;n-1,m-1) \,\sspp\, (a\sspp d\,m)\, S2fac(d,a;n-1,m), \ \ {\rm for}\ n\sspgeq 1\,,\ m\sspeq 0,\,1,\,...,\,n\,,  
\Eeq
with $S2fac(d,a;n,-1)\sspeq 0$, $ S2fac(d,a;n,m)\sspeq 0$ for $n\sspkl m$ and $S2fac(d,a;0,0)\sspeq 1$.\psn
The \egf of the row polynomials of ${\bf S2fac}[d,a]$ is \dstyle{\frac{e^{a\,t}}{1\sspm (e^{d\,t\sspm 1})}}. This is seen after interchanging the two sums and using the \egf of columns of ${\bf S2}[d,a]$.\psn  
For the instances $[d,a]\speq [1,0],\,[2,1],\,[3,1],\, [3,2],\,[4,1],\, ]4,3]$ see \seqnum{A131689}, \seqnum{A145901},\, \seqnum{A284861},\, \seqnum{A225472},\, \seqnum{A285066},\, \seqnum{A225473}, respectively. \pbn
{\bf 2. Proof of eq.\,${\bf (13)}$, i.e., eq.\,${\bf (12)}$}\psn
The {\sl Sheffer} structure eq.\,$(12)$ means that the column \egf of the ${\bf S2}[d,a]$ number triangle satisfies eq.\,$(13)$. The column \egf $ES2Col(d,a;t,m)$ is here named $E(t,m)$ for simplicity.
The lower summation index in brackets can be used instead of the given one because of the triangle structure of ${\bf S2}[d,a]$.  The recurrence is used in the first step.
\Beqarray
E(t,m) &\sspdef& \sum_{n=m(0)}^{\infty}\,S2(d,a;n,m)\frac{t^n}{n!} \nonumber\\
 &\sspeq& d\,\sum_{n=0(1)}^{\infty}\, S2(d,a;n-1,m-1)\frac{t^n}{n!} \sspp (a\sspp d\,m)\,\sum_{n=0(1)}^{\infty}\,S2(d,a;n-1,m)\,\frac{t^n}{n!} \nonumber\\
&\sspeq& \int\,dt\,\left(d\,\sum_{n=1}^{\infty}\, S2(d,a;n-1,m-1)\frac{x^{n-1}}{(n-1)!} \sspp (a\sspp d\,m)\,\sum_{n=1}^{\infty}\,S2(d,a;n-1,m)\,\frac{t^{n-1}}{(n-1)!} \right)\ \nonumber \\
&\sspeq& \int\,dt\,\left(d\,E(t,m-1) \sspp (a\sspp d\,m)\,E(t,m)\right)\ .
\Eeqarray
Differentiating both sides of the final equation produces a recurrence for $E(t,m)$:
\Beq
\left(\frac{d\ }{dt} \sspm (a\sspp d\,m)\right)\,E(t,m)\sspeq d\,E(t,m-1)\,, 
\Eeq
with the input $E(t,0)\sspeq e^{a\,t}$ because $S2(d,a;n,0)\sspeq a^n$ from the recurrence.\psn
The solution of this differential difference equation, satisfying the input, is
\Beq\
E(t,m) \sspequiv ES2Col(d,a;t,m)\sspeq e^{a\,t}\,\frac{\left(e^{d\,t}\sspm 1\right)^m}{m!}\,, 
\Eeq
which is eq.\,($13$).\pbn
{\bf 3. Proof of eq.\,$\bf (14)$}\psn
The \ogf of eq.\,$(14)$ $GS2Col(d,a;x,m) \sspequiv  G(x,m)$ for short in this proof, is shown to lead to the \egf $ES2Col(d,a;t,m)\sspequiv E(t,m)$ of eq.\,$(13)$ \via\ the inverse {\sl Laplace} transform, given in eq.\,$(4)$. \psn
{\bf Lemma 2:} 
\Beq
\prod_{j=0}^m\,\frac{1}{x\sspm (a\sspp d\,j)}\sspeq \left(\frac{-1}{d}\right)^m\,\frac{1}{m!}\,\sum_{j=0}^m\,(-1)^j\,{\binomial{m}{j}}\,\frac{1}{x\sspm (a\sspp d\,j)}\, ,\ \ {\rm for}\ \ m\sspin \mathbb N_0\, . 
\Eeq
{\bf Proof}: This is a standard partial fraction decomposition for the rational function \dstyle{\frac{1}{P(x)}} with $P(x)$ a polynomial of degree $m+1$ with the simple roots $\alpha_j\sspeq a\sspp d\,j$, $j\sspeq 0,\,1,\,...,\,m $. \dstyle{\frac{1}{P(x)} \sspeq \sum_{j=0}^m\,\frac{a_j}{x \sspm (a\sspp d\,j)}}. Here \dstyle{a_j\sspeq \frac{1}{P^{\prime}(\alpha_j)}\sspeq \prod_{k=0,\neq j}^m\, \frac{1}{\alpha_j\sspm \alpha_k}\sspeq\frac{1}{d^m}\, \prod_{k=0,\neq j}^m\,\frac{1}{j-k}\sspeq \frac{1}{d^m}\,\frac{1}{j!\,(-1)^{m-j}\,(m-j)!}\sspeq \frac{1}{d^m}\,\frac{1}{m!}\,(-1)^{m-j}\,{\binomial{m}{j}}}.\hskip .7cm $\square$\psn
Now, due to the linearity of ${\cal L}^{-1}$, and the transform \dstyle{ {\cal L}^{-1}\left[\frac{1}{p\sspm\alpha} \right]\sspeq e^{\alpha\,t}}, one finds after application of {\sl Lemma 2}, and with the help of the binomial theorem:
\Beqarray
E(t,m)&\sspeq& {\cal L}^{-1}\left[\frac{1}{p}\,G\left(\frac{1}{p}\,,m\right)\right]\sspeq {\cal L}^{-1}\left[\frac{d^m}{\prod_{j=0}^m\,(p\sspm (a\sspp d\, j))} \right] \sspeq \frac{(-1)^m}{m!}\,\sum_{j=0}^m\,(-1)^j\,{\binomial{m}{j}}\,{\cal L}^{-1}\left[ \frac{1}{p\sspm (a\sspp d\,j)}\right]    \nonumber \\
&\sspeq& \frac{(-1)^m}{m!}\,e^{a\,t}\,\sum_{j=0}^m\,{\binomial{m}{j}}\, (-e^{d\,t})^j\sspeq \frac{e^{a\,t}}{m!}\,(1\sspm e^{d\,t})^m \,(-1)^m\, .
\Eeqarray
which is indeed the \egf given in eq.\,$(13)$.\pbn
{\bf 4. Proof of eq.\,$\bf (16)$ with eq.\,$\bf (17)$}\psn
{\bf Lemma 3: Sheffer transform of a sequence}\psn
If the \egf  of sequence $\{b_n\}_{n=0}^{\infty}$ is ${\cal B }(t)$, the \egf of sequence $\{a_n\}_{n=0}^{\infty}$ is ${\cal {A}}(t)$, and the {\sl Sheffer} transform of $\{a_n\}$ is $b_n \sspeq \sum_{m=0}^n\,S(n,\, m)\,a_n$, with $S$ {\sl Sheffer} of type $S\sspeq (g(t),\ f(t))$ then 
\Beq
{\cal B}(t)\sspeq g(t)\,{\cal A}(f(t))\ .
\Eeq
The {\bf proof} uses an exchange of the two summations and the \egf of column $m$ of $S$, \ie \dstyle{g(t)\,\frac{f(t)^m}{m!}}.\psn
{\bf  Corollary 1}: \psn
The row polynomials $PS(n,\,x)$ of a {\sl Sheffer} matrix ${\bf S}\sspeq (g(t),\,f(t))$ have \egf \dstyle{EPS(t,\,x)\sspeq g(t)\,e^{x\,f(t)}}. \pn
This is also called the \egf of the triangle $\bf S$.
\pbn 
That the  $fallfac(d,a;x,m)$ definition in \Eq{17} can be rewritten in terms of the usual falling factorial $x^{\underline{n}}$ is trivial.\psn
\vfill
\eject
\noin
{\bf Lemma 4: E.g.f. of fallfac[d,a]}\psn
The \egf of $fallfac(d,a;x,m)$ (see eq. $(17)$) is
\Beq
F(d,a;x,t)\sspdef 1\sspp \sum_{m=1}^{\infty}\, fallfac(d,a;x,m)\,\frac{t^m}{m!}\sspeq (1\sspp d\,t)^{\frac{x-a}{d}}\, .  
\Eeq
{\bf Proof}: With the binomial theorem and the rewritten form of $fallfac[d,a]$ in terms of the ordinary falling factorial this is trivial.\psn
Now eq.\,$(16)$ is a {\sl Sheffer} transform of the sequence $\{fallfac(d,a;x,m)\}_{m=0}^{\infty}$, therefore, with {\sl Lemmata} $3$ and $4$
the \egf of the \rhs of eq.\,$(16)$ is, with the \egf of ${\bf \widehat{S2}}$ given in connection with eq.\,$(15)$,
\Beq
e^{a\,t}\,F\left(d,a;t,\frac{1}{d}\,(e^{d\,t}\sspm 1)\right)\sspeq e^{a\,t}\, \left(1\sspp (e^{d\,t}\sspm 1) \right)^{\frac{x-a}{d}}\sspeq e^{a\,t}\,e^{d\,t\,\frac{x-a}{d}}\sspeq e^{t\,x},
\Eeq
which is the \egf of the sequence $\{x^n\}_{n=0}^{\infty}$ of the \lhs.\pbn
{\bf 5. Meixner recurrence and recurrence for Sheffer row polynomials eq.\,$\bf (19)$}\psn
{\bf a)} Monic row polynomials $s(n,\,x)$ of the {\sl Sheffer} type $(g(x),\, f(x))$ satisfy the {\sl Meixner} \cite{Meixner}, p. 9, eqs.\,$(4.1)$ and $(4.2)$, recurrence ($f^{[-1]}$ denotes the compositional inverse of $f$)
\Beq
f^{[-1]}({\bf d}_x)\, s(n,\,x) \sspeq n\, s(n-1,\,x),\ \ {\rm with\ input}\ \ s(0,\,x)\sspeq 1.
\Eeq
For the proof see the original reference.\psn
For $P\widehat{S2}(d,a;n,x)\sspeq \sum_{m=0}^n\,\widehat{S2}(d,a;n,m)\,x^m$ with \dstyle{f^{[-1]}(y)\sspeq \frac{1}{d}\,log(1\sspp\,y)}
one has 
\Beq
\frac{1}{d}\,\sum_{k=1}^n\, \frac{(-1)^{k+1}}{k}\, \frac{d^k\,\, }{dx^k}\,P\widehat{S2}(d,a;n,x)\sspeq n\,P\widehat{S2}(d,a;n-1,x),\ \ {\rm with\ input}\ \ P\widehat{S2}(d,a;0,x)\sspeq 1\ .
\Eeq
{\bf b)} The standard recurrence for {\sl Sheffer} row polynomials $S(n,\,x)$ (not necessarily monic ones) is given in {\sl Roman} \cite{Roman}, p. 50, Corollary 3.7.2, which in our notation is\psn
{\bf Lemma 5:} (the differentiation is with respect to $t$)
\Beq
S(n,\, x)\sspeq \left.\left[x\sspp \left(log(g(f^{[-1]}(t)))\right)^{\prime}\right]\,\frac{1}{f^{[-1]}(t)^{\prime}}\,\right |_{t\sspeq {\bf d}_x}\, S(n-1,\,x), \ \ {\rm for}\ n\sspin \mathbb N\,,
\Eeq
and input $S(0,\,x)\sspeq 1 $.\psn
For $PS2(d,a;n,x)$ of eq.\,$(19)$ $f^{[-1]}(t)\sspeq \frac{1}{d}\,log(1\sspp t )$, $ f^{[-1]}(t)^{\prime}\sspeq \frac{1}{d}\,\frac{1}{1\sspp t}$, $g(t)\sspeq e^{a\,t}$ and 
$\left(log(g(f^{[-1]}(t)))\right)^{\prime}\sspeq \frac{a}{d}\,\frac{1}{1\sspp t}$, leading to \dstyle{PS2(d,a;n,x)\sspeq \left. [x\,d\,(1\sspp t)\sspp a]\right |_{t\sspeq {\bf d}_x}\,PS2(d,a;n-1,x)}, which his eq.\,$(19)$,\pbn
{\bf 6. Proof of eqs.\,$\bf (23)$ to $\bf (25)$}\psn
Multiplication of eq.\,$(23)$ with $(1\sspm x)^{n+1}$ and the binomial formula gives
\Beqarray
\sum_{i=0}^n\,a_i^{(n)}\, x^i&\sspeq& \sum_{j=0}^n\, b_j^{(n)}\,x^j\,(1\sspm x)^{n-j}\sspeq \sum_{j=0}^n\,b_j^{(n)}\,\sum_{k=0}^{n-j}\,(-1)^k\,{\binomial{n-j}{k}}\, x^{k+j}\nonumber \\
&\sspeq& \sum_{i=0}^n\,x^i\left(\sum_{j=0}\,b_j^{(n)}\,(-1)^{i-j}\,{\binomial{n-j}{i-j}} \right)\, ,
\Eeqarray
where in the last step a new summation index $i\sspeq k+j$ has been used instead of $k$, and the upper summation index $j$ is determined by the binomial as $\min(n,\,i)\sspeq i$, because $0\sspleq i\sspleq n$. Comparing the coefficients of the powers $x^i$, for $i\sspeq 1,\,2,\, ...,\,n$, leads then to eq.\,$(24)$ for $a_i^{(n)}$.\psn
The inverse relation, eq.\,$(25)$, uses the following binomial identity (see \cite{GKP}, p. 169. eq.\,$(5.24)$ with $k\sspto p$, $m\sspeq s\sspeq 0$, $l\sspto n-k$, $n\sspto n \sspm j$) 
\Beq
\sum_{p\sspgeq 0}\,(-1)^p\,{\binomial{n-k}{p}}\,{\binomial{p}{n-j}}\sspeq (-1)^{n-k}\,{\binomial{0}{k-j}}\sspeq (-1)^{n-k}\,\delta_{k,j}\,,
\Eeq
with the {\sl Kronecker} symbol $\delta_{j,k}\sspeq 1$ if $j=k$, and $0$ otherwise.\psn
In the \rhs of eq.\,$(25)$, with $a_i^{(n)}$ from eq.\,$(24)$ inserted, the two finite sums are interchanged, a new summation index $p\sspeq n-i$ is used instead of $i$, and finally the preceding binomial identity is employed.
\Beqarray
\sum_{i=0}^j\,{\binomial{n-i}{j-i}}\,a_i^{(n)}&\sspeq& \sum_{i=0}^j\,{\binomial{n-i}{j-i}}\,\sum_{k=0}^i\,(-1)^{i-k}\,b_k^{(n)}\,{\binomial{n-k}{i-k}} \sspeq \sum_{k=0}^j\,(-1)^{i-k}\,b_k^{(n)}\,\left(\sum_{i=k}^j\,{\binomial{n-i}{n-j}}\, {\binomial{n-k}{n-i}}\right)\nonumber \\
&\sspeq& \sum_{k=0}^j\,(-1)^{n-k}\,b_k^{(n)}\,\left(\sum_{p=n-j}^{n-k}\,(-1)^p\,{\binomial{n-k}{p}}\,{\binomial{p}{n-j}}\right)\sspeq b_j^{(n)}\, . 
\Eeqarray
{\bf 7. Proof of eqs.\,$\bf (26)$ to $\bf (28)$}\psn
This follows by using in eq.\,$(21)$ eqs.\, $(23)$ with $b_j^{(n)}\sspeq S2(d.a;n,j)\,j!$ and eq.\,$(24)$. Then $a_k^{(n)}\sspeq rEu(d,a;n,k)$ is obtained as given in eq.\,$(28)$.\pbn
{\bf 8. Proof of eq.\,$\bf (29)$}\psn
This is eq.\,$(25)$ (replacing $i\sspto k$) with $b_j^{(n)}$ and $a_k^{(n)}$ as given in the previous proof.\pbn
{\bf 9. Proof of eq.\,$\bf (30)$}\psn
This uses the binomial identity (see \cite{GKP}, p. 169, eq. $(5.26)$ with $k\sspto j$, $q\sspto 0$, $m\sspto n-k$, $l\sspto n$ and $n\sspto l$)\psn
\Beq
\sum_{j=l}^k\,{\binomial{n-j}{n-k}}\, {\binomial{j}{l}}\sspeq  {\binomial{n+1}{n-k+l+1}}\sspeq  {\binomial{n+1}{k-l}}\, .
\Eeq
Insertion of eq.\,$(11)$ into eq.\,$(28)$ followed by an interchange of the two sums and application of the binomial identity leads to
\Beqarray
rEu(d,a;n,k)&\sspeq& \sum_{j=0}^k\,(-1)^{k-j}\,{\binomial{n-j}{k-j}}\, \sum_{l=0}^j\,(-1)^{j-l}\,{\binomial{j}{l}}\,(a\sspp d\,l)^n\sspeq \sum_{l=0}^k\,(-1)^{k-l}\, (a\sspp d\,l)^n\,\sum_{j=l}^k\,{\binomial{n-j}{n-k}}\, {\binomial{j}{l}} \nonumber \\
&\sspeq& \sum_{l=0}^k\,(-1)^{k-l}\,(a\sspp d\,l)^n\,{\binomial{n+1}{k-l}}\sspeq \sum_{j=0}^k\, (-1)^{k-j}\,{\binomial{n+1}{k-j}}\,(a\sspp d\,j)^n\, .
\Eeqarray
For the ${\bf rEu}[d,a]$ triangles with $[d,a]\sspeq [1,0],\, [2,1],\, [3,2],\,[4,3]$ see \seqnum{A123125}, \seqnum{A060187}, \seqnum{A225117}, \seqnum{A225118}, respectively. The case ${\bf rEu}[3,1]$ is the row reversed version of ${\bf rEu}[3,2]$, and ${\bf rEu}[4,1]$ is the row reversed version of ${\bf rEu}[4,3]$, In general this row reversion relation holds  between ${\bf rEu}[d,d-a]$ and ${\bf rEu}[d,a]$, for \dstyle{a\sspeq 1,\,...,\,\floor{\frac{d}{2}}} for $\gcd(d-a,a)\sspeq 1$.\pbn
{\bf 10. Proof of eq.\,$\bf (31)$}\psn
Like eq.\,$(6)$ one has for $GP(d,a;n,x)$ of eq.\,$(20)$
\Beq
GP(d,a;n,x)\speq \sum_{k=0}^n\,{\binomial{n}{k}}\,a^{n-k}\,d^k\,GP(k,x)\,,
\Eeq
with $GP(k,x)\sspeq GP(1,0;k,x)$. However, this does not lead immediately to the desired formula for $rEu(d,a;n,k)$ in terms of the usual {\sl Euler}ian numbers $rEu(n,\,k)$ claimed in eq.\,$(31)$. The proof is done by inserting ${\bf S2}[d,a]$ from eq.\,$(18)$ into eq.\,$(28)$, then replacing the usual ${\bf S2}$ by ${\bf rEu}$ \via\ eq.\,$(29)$ for $[d,a]\sspeq [1,0]$. Here two binomial identities are needed. The first is given in \cite{GKP}, p. 169, eq.\,$(5.25)$ (with $k\sspto j$, $s\sspto m-i$, $n\sspto i$, $l\sspto n$, and $m\sspto n-k$). Here one needs $n-k\sspgeq 0$ and the upper summation index which would be $n$ can be replaced by $k$ because for $j\sspeq k=1,\,...,\,n$ the first binomial vanishes because the upper non-negative number is then smaller than the lower one.  
\Beq
\sum_{j=i}^k\,(-1)^j\, {\binomial{n-j}{n-k}}\,{\binomial{m-i}{j-i}}\sspeq (-1)^k\, {\binomial{m-n+k-i-1}{k-i}}\sspeq (-1)^k\, {\binomial{-(n-m-(k-i)+1)}{k-i}}\sspeq (-1)^{i}\,{\binomial{n-m}{k-i}}\ .
\Eeq
In the last step another identity, given in \cite{GKP}, p. 164, \Eq{5.14}, has been used.\psn
\Beq
{\binomial{-r}{p}} \sspeq (-1)^p\,{\binomial{p-1+r}{p}}\ .
\Eeq
Now  
\Beqarray
rEu(d,a;n,k)&\sspeq& \sum_{j=0}^k\,(-1)^{k-j}\,{\binomial{n-j}{k-j}}\,S2(d,a;n,j)\,j!\sspeq \sum_{j=0}^k\,(-1)^{k-j}\,{\binomial{n-j}{k-j}}\,\sum_{m=0}^n\,{\binomial{n}{m}}\,a^{n-m}\,d^m\,S2(m,\,j)\,j!\nonumber \\
&\sspeq& \sum_{j=0}^k\,(-1)^{k-j}\,{\binomial{n-j}{k-j}}\,\sum_{m=0}^n\,{\binomial{n}{m}}\,a^{n-m}\,d^m\,\sum_{i=0}^j\,{\binomial{m-i}{j-i}}\,rEu(m,\,i)\nonumber \\
&\sspeq& \sum_{m=0}^n\,{\binomial{n}{m}}\,a^{n-m}\,d^m\,\sum_{i=0}^k\,(-1)^k\,rEu(m,\,i)\,\sum_{j=0}^k\,(-1)^j\, {\binomial{n-j}{n-k}}\,{\binomial{m-i}{j-i}}\nonumber\\
&\sspeq& \sum_{m=0}^n\,{\binomial{n}{m}}\,a^{n-m}\,d^m\,\sum_{i=0}^k\,(-1)^k\,rEu(m,\,i)\,(-1)^i\, {\binomial{n-m}{k-i}}\, ,
\Eeqarray
which is eq.\,$(31)$ with summation index $p\sspto i$.\pbn
{\bf 11. Proof of eq.\,$\bf (32)$}\psn
We show that eq.\,$(30)$ satisfies the three term recurrence eq.\,$(31)$. The \lhs of the recurrence is 
\Beq
rEu(d,a;n,m)\sspeq \sum_{j=0}^m\,(-1)^{m-j}\, {\binomial{n+1}{m-j}}\,(a\sspp d\,j)^n\, .
\Eeq
The \rhs of the recurrence is
\Beqarray
&&(d\,(n-m)\sspp (d\sspm a))\,\sum_{j=0}^{m-1}\,(-1)^{m-j-1}\,{\binomial{n}{m-1-j}}\,(a\sspp d\,j)^{n-1}\nonumber\\
&&\sspp (a\sspp d\,m)\,\sum_{j=0}^m\,(-1)^{m-j}\,{\binomial{n}{m-j}}\,(a\sspp d\,j)^{n-1}\nonumber\\ 
&\sspeq& \sum_{j=0}^m\,(-1)^{m-j}\,{\binomial{n+1}{m-j}}\,(a\sspp d\,j)^{n-1} \left[-(d\,(n-m)\sspp (d\sspm a))\frac{{\binomial{n}{m-1-j}}}{{\binomial{n+1}{m-j}}}\sspp  (a\sspp d\,m)\,\frac{{\binomial{n}{m-j}}}{{\binomial{n+1}{m-j}}}\right]\, .  
\Eeqarray
In the first sum the upper index $m-1$ has been extended to $m$ because the extra term vanishes due to the binomial. The terms in the bracket are shown to become $a\sspp d\,j$ as follows.\psn
\Beqarray
[...]&\sspeq& -(d\,(n-m)\sspp (d\sspm a))\,\frac{n!\,(m-j)!}{(m-1-j)!\,(n+1)!}\sspp (a\sspp d\,m)\, \frac{n!(n+1-m+j)!}{(n-m+j)!\,(n+1)!} \nonumber \\
&\sspeq& \frac{1}{n+1}\,\left(-(d\,(n-m)\sspp (d\sspm a))\,(m-j)\sspp (a\sspp d\,m)\,(n-m+j+1)\right)\nonumber \\
&\sspeq& \frac{1}{n+1}\, \left( d\,j\,(n+1)\sspp a\,(n+1)\right)\sspeq a+d\, j\,.
\Eeqarray
{\bf 12. Proof of eq.\,$\bf (34)$}\psn
From eq.\,(21) and the row polynomial definition of $PS2fac$ we have
\Beq
GP(d,a;n,x)\sspeq \frac{1}{1-x}\,\sum_{k=0}^n\,S2(d,a;n,k)\,k!\,\left(\frac{x}{1-x}\right)^k\sspeq \spfed \frac{1}{1-x}\,PS2fac\left(d,a;n,\frac{x}{1-x}\right)\ .
\Eeq
Therefore, from eq.\,$(26)$, 
\Beq
\frac{1}{1-x}\,PS2fac\left(d,a;n,\frac{x}{1-x}\right)\sspeq \frac{1}{(1-x)^{n+1}}\,PrEu(d,a;n,x)\, ,
\Eeq
which is \Eq{34} for $x\sspneq 1$. But \Eq{34} holds also for $x\sspeq 1$, with the row sums $PrEu(d,a;n,1)\sspeq S2fac(d,a;n,n)$. From \Eq{14} one has $S2fac(d,a;n,n) \sspeq [x^n]\left( n!\,GS2Col(d,a;x,n) \right)\sspeq d^n\,n!$ with \egf \dstyle{\frac{1}{1\sspm d\,t}}, independently of $a$. These sequences are, for $d=1,\,2,\,...,\,5$, \seqnum{A000142}, \seqnum{A000165}, \seqnum{A032031}, \seqnum{A047053}, \seqnum{A052562}\,.
\pbn
{\bf 13. Proof of eq.\,$\bf (35)$}\psn
With eq.\,$(34)$ and the {\sl Sheffer} structure of ${\bf S2}[d,a]$ (eq.\, $(13)$ for the column \egf, here in the variable $(1-x)\,z$)
\Beqarray
EPrEu(d,a;t,x)&\sspdef& \sum_{n=0}^{\infty}\,\frac{t^n}{n!}\,PrEu(d,a;n,x)\sspeq \sum_{n=0}^{\infty}\,\frac{(t\,(1-x))^n}{n!}\,\sum_{m=0}^n\,S2(d,a;n,m)\,m!\left(\frac{x}{1-x}\right)^m \nonumber \\
&\sspeq& \sum_{m=0}^{\infty}\,\left(\frac{x}{1-x}\right)^m\, \sum_{n=m(0)}^{\infty}\, \frac{(t\,(1-x))^n}{n!}\,S2(d,a;n,m)\,m!\nonumber \\
&\sspeq& e^{a\,(1-x)\,t}\,\sum_{m=0}^{\infty}\,\left(\frac{x}{1-x}\right)^m \frac{(e^{d\,(1-x)\,t}\sspm 1)^m}{m!}\, m!\nonumber \\
&\sspeq&  e^{a\,(1-x)\,t}\,\frac{1}{1\sspm \frac{x}{1-x}\,(e^{d\,(1-x)\,t}\sspm 1)}\sspeq \frac{(1-x)\,e^{a\,(1-x)\,t}}{1\sspm x\,e^{d\,(1-x)\,t}}\, .
\Eeqarray
The limit $x\sspto 1$ \via\ {\sl l'H\^opital}'s rule leads to the \egf \dstyle{\frac{1}{1\sspm d\,t}} for the row sums of ${\bf rEu}[d,a]$, as found also in the preceding proof. 
\pbn
{\bf 14. Proof of eq.\,$\bf (39)$}\psn
One obtains $b_j^{(n+1)}$ of eq.\,$(37)$ from eq.\,$(25)$ (with $n\sspto n+1$) with $a_i^{(n+1)}$ given in eq.\,$(38)$. We omit the $(d,a)$ labels here. Remember that $a_{n+1}^{(n+1)}\sspeq 0$ (but $b_{n+1}^{(n+1)}$ does not vanish).
\Beq
b_j^{(n+1)}\sspeq \sum_{i=0}^j\,{\binomial{n+1-i}{j-i}}\,a_i^{(n+1)}\sspeq \sum_{i=0}^j\,{\binomial{n+1-i}{j-i}}\,\sum_{p=0}^i\,(-1)^{i-p}\,{\binomial{n+1}{i-p}}\,(a\sspp d\,p)^n\, .
\Eeq
It is convenient to consider the case  $j=0$ separately.
\Beq
 b_0^{(n+1)}\sspeq a_0^{(n+1)}\sspeq a^n,\  \ n\sspin \mathbb N_0\, .
\Eeq
For $j\sspeq 1,\,2,\,...,\,n+1$ we have after exchange of the sums
\Beq
b_j^{(n+1)}\sspeq \sum_{p=0}^j\,(-1)^p\, (a\sspp d\,p)^n\,\sum_{i=p\,(0)}^{j}\,(-1)^i\,{\binomial{n+1-i}{n+1-j}}\,{\binomial{n+1}{i-p}}\,. 
\Eeq
Because of the second binomial one could start the sum with $i=0$. Now the binomial identity used already above in the first step of eq.\,$(95)$, is employed with $j\sspto i,\, i\sspto p,\,n\sspto n+1,\, k\sspto j,\, m\sspto p+n+1$. In the first binomial the lower number is non-negative and the original upper sum index $n+1$ can be replaced by $j$ because for $i\sspeq j+1,\,...,\,n+1$ the upper non-negative number in the first binomial becomes smaller then the lower one.
\Beq
b_j^{(n+1)}\sspeq \sum_{p=0\,(1)}^j\,(-1)^p\, (a\sspp d\,p)^n\,(-1)^j\,{\binomial{j-1}{j-p}}\, . 
\Eeq
Because $j\sspgeq 1$ one can use the symmetry of the binomials (this is why we have separated the $j=0$ case).
\Beq
b_j^{(n+1)}\sspeq \sum_{p=0}^j\,(-1)^{j-p}\,{\binomial{j-1}{p-1}}\,(a\sspp d\,p)^n\, . 
\Eeq
The identification with the $\Sigma S2[d,a]$ as claimed in eq.\,$(39)$ is achieved by using the {\sl Pascal} recurrence (see \seqnum{A007318}).
\Beqarray
&&\sum_{p=0}^j\,(-1)^{j-p}\,{\binomial{j-1}{p-1}}\,(a\sspp d\,p)^n\sspeq - \sum_{p=0}^j\,(-1)^{j-p}\,{\binomial{j-1}{p}}\,(a\sspp d\,p)^n \sspp \sum_{p=0}^j\,(-1)^{j-p}\,{\binomial{j}{p}}\,(a\sspp d\,p)^n \nonumber\\
&&\sspeq +\, S2(d,a;n,j-1)\,(j-1)!\sspp S2(d,a;n,j)\,j!\, ,
\Eeqarray
where eq\,$(11)$ has been used. Now the $j=0$ result $a^n$ is also covered because  $S2(d,a;n,0)\sspeq a^n$ and  $S2(d,a;n,-1)*(-1)!$ is taken as vanishing ($S2(d,a;n,-1)=0$ from the recurrence eq.\,$(10)$).
\vfill
\eject
\noin
{\bf 15. Proof of eq.\,$\bf (40)$}\psn
Putting things together, the \egf of $\{PS(d,a;n,m)\}_{m=0}^{\infty}$ (see \Eq{1}) becomes \via\ inverse {\sl Laplace} transform of \Eq{2}
\Beq
EPS(d,a;n,t)\sspdef \sum_{m=0}^{\infty}\,PS(d,a;n,m)\,\frac{t^m}{m!} \sspeq {\cal L}^{-1}\left [\frac{1}{p}\,GPS\left(d,a;n,\frac{1}{p}\right) \right ]\,,
\Eeq
and from eq.\,$(36)$ with eqs.\,$(27)$,\,$(38)$ and $(37)$
\Beq
\frac{1}{p}\,GPS\left(d,a;n,\frac{1}{p}\right)\sspeq \sum_{j=0}^{n+1}\,b^{(n+1)}_j(d,a)\, \frac{1}{p}\,\frac{1}{p^j}\,\frac{1}{\left(1\sspm \frac{1}{p}\right)^{j+1}}\sspeq \sum_{j=0}^{n+1}\,b^{(n+1)}_j(d,a)\,\frac{1}{(p\sspm 1)^{j+1}}\, .
\Eeq
Thus, by linearity of ${\cal L}^{-1}$, and the formula before eq.\,$(22)$, one obtains
\Beq
EPS(d,a;n,t)\sspeq e^t\,\sum_{j=0}^{n+1}\,b_j^{(n+1)}(d,a)\,\frac{t^j}{j!}\,,
\Eeq
which become finally eq.\,$(40)$ after insertion of $b_j^{(n+1)}(d,a)$ from eq.\,$(39)$.
\pbn
{\bf B) Proofs of section 1\,B}\psn
{\bf 1. Proof of eq.\,$\bf (47)$}\psn
Insertion of eq.\,$(18)$ into eq.\,$(46)$ leads after interchange of the two finite sums to
\Beq
B(d,a;n)\sspeq \sum_{k=0}^n\,{\binomial{n}{k}}\,a^{n-k}\,d^k\,\sum_{m=0}^{n\ (k)}\,(-1)^m\,\frac{1}{m+1}\,S2(k,\,m)\,m!\,,
\Eeq
where in the second sum the upper index can be taken as $k$ instead of $n$ because $S2(k,\,m)\sspeq 0$ for $m\sspgr k$, and $k\sspleq n$ from the first sum. Then the second sum is equal to $B(k)$ by eq.\,$(46)$ for $[d,a]\sspeq [1,0]$. This is then eq.\,$(47)$.\psn
For the $B[d,a]$ numbers for $[d,a] = [1,0],\, [2,1],\,[3,1],\,[4,1],\,[5,1],\,[5,2]$
see \seqnum{A027641}/\seqnum{A027642}, \seqnum{A157779}/\seqnum{A141459}, \seqnum{A157799}/\seqnum{A285068},
\seqnum{A157817}/\seqnum{A141459}, \seqnum{A157866}/\seqnum{A288872}, \seqnum{A157833}/\seqnum{A288872}, respectively.
\pn
$B[3,2](n)\sspeq (-1)^n\,B[3,1](n)$, $B[4,3](n)\sspeq (-1)^n\,B[4,1](n)$, $B[5,3](n)\sspeq (-1)^n\,B[5,2](n)$,
and $B[5,4](n)\sspeq (-1)^n\,B[5,1](n)$. \pbn
{\bf 2. Proof of eq.\,$\bf (48)$}\psn
The \egf of $\{B(d,a;n)\}_{n=0}^\infty$ is obtained from the defining eq.\,$(46)$ recognizing, after exchange of the sums, the \egf  $ES2Col(d,a;t,m)$ of eq.\,$(13)$:\psn
\Beqarray
EB(d,a;t) &\sspdef& \sum_{n=0}^{\infty}\, B(d,a;n)\,\frac{t^n}{n!}\sspeq \sum_{n=0}^{\infty}\, \frac{t^n}{n!}\, \sum_{m=0}^n\, (-1)^m\, \frac{m!}{m+1}\, S2(d,a;n,m)\nonumber\\
&\sspeq& \sum_{m=0}^{\infty}\, (-1)^m\, \frac{m!}{m+1}\, \sum_{n=m}^{\infty}\, \frac{t^n}{n!}\, S2(d,a;n,m)\sspeq  \sum_{m=0}^{\infty}\, (-1)^m\, \frac{1}{m+1}\, e^{a\,t}\,(e^{d\,t}\sspm 1)^m \nonumber \\
&\sspeq& e^{a\,t}\, \frac{1}{y}\,\sum_{m=0}^{\infty} \frac{y^{m+1}}{m+1}\ \ ({\rm with} \ y\sspdef 1\sspm e^{d\,t})
\sspeq e^{a\,t}\, \frac{1}{y}\,\int\,dy \sum_{m=0}^{\infty}\, y^m \sspeq e^{a\,t}\, \frac{1}{y}\,\int\,\frac{dy}{1-y}\nonumber \\
&\sspeq& e^{a\,t}\, \frac{1}{-y}\,log(1-y)\sspeq \frac{d\,t\,e^{a\,t}}{e^{d\,t}\sspm 1}\,.
\Eeqarray
{\bf 3. Proof of eq.\,$\bf (50)$}\psn
This follows from inserting into  \Eq{49} (with $n-m\sspeq p$) eq,\,$(47)$, using the binomial identity (see \cite{GKP}, p. 174. Table 174. the trinomial revision formula)\dstyle{{\binomial{n}{p}}\,{\binomial{p}{m}}\sspeq {\binomial{n}{m}}\, {\binomial{n-m}{p-m}}} and interchange of the sums. Then the binomial formula is used.\psn
\Beqarray
B(d,a;n,x) &\sspeq & \sum_{p=0}^n\,{\binomial{n}{p}}\,B(d,a;p)\,x^{n-p}\sspeq \sum_{p=0}^n\,{\binomial{n}{p}}\,x^{n-p}\, \sum_{m=0}^p\,{\binomial{p}{m}}\,a^{p-m}\,d^m\,B(m)\nonumber \\
 &\sspeq & \sum_{p=0}^n\,\sum_{m=0}^p\, {\binomial{n}{m}}\, {\binomial{n-m}{p-m}}\, x^{n-p}\,a^{p-m}\,d^m\,B(m)\sspeq \sum_{m=0}^n\,{\binomial{n}{m}}\,a^{-m}\,d^m\,B(m)\,\sum_{p=m}^n\,{\binomial{n-m}{p-m}}\, x^{n-p}\,a^p \nonumber \\
 &\sspeq & \sum_{m=0}^n\,{\binomial{n}{m}}\,a^{-m}\,d^m\,B(m)\,\sum_{p=0}^{n-m}\,{\binomial{n-m}{p}}\, x^{n-(p+m)}\,a^{p+m}\nonumber \\
&\sspeq & \sum_{m=0}^n\,{\binomial{n}{m}}\,d^m\,B(m)\,\sum_{p=0}^{n-m}\,{\binomial{n-m}{p}}\, x^{(n-m)-p}\,a^p\sspeq \sum_{m=0}^n\,{\binomial{n}{m}}\,d^m\,B(m)\,(a\sspp x)^{n-m}\, .
\Eeqarray
{\bf 4. Proof of eq.\,$\bf (51)$}\psn
Eq.\,$(49)$ is an exponential (also called binomial) convolution of the sequences $\{B(d,a;n)\}_{n=0}^{\infty}$ and $\{x^n\}_{n=0}^{\infty}$, hence the product of their {\sl e.g.f.}s is with \Eq{48} \dstyle{EB(d,a;t)\,e^{x\,t}\sspeq \frac{d\,t\,e^{(a+x)\,t}}{e^{d\,t}\sspm 1}}.\pn
Alternatively one can take the exponential convolution of eq.\,$(50)$ of $\{d^m\,B(m)\}_{m=0}^{\infty}$ with \egf \dstyle{\frac{d\,t}{e^{d\,t}\sspm 1}} and $\{(a\sspp x)^n\}_{n=0}^{\infty}$ with \egf $e^{(a+x)\,t}$, and their product is eq.\,$(51)$.\pbn
{\bf 5. Proof of eq.\,$\bf (57)$}\psn
This proof will be rather lengthy. It will need the following simple {\it Lemma}.\psn
{\bf Lemma 6:} If the \egf of the sequence $\{C_n\}_{n=0}^{\infty}$ is ${\cal C}(t)$ then the \egf of the sequence \dstyle{ \left\{\frac{C_{n+1}}{n+1}\right\}_{n=0}^{\infty}} is \dstyle{\frac{1}{t}\,( {\cal C}(t)\sspm C_0)}\ .
\psn
{\bf Proof:} \pn
\Beq
\sum_{n=0}^{\infty}\,  \frac{C_{n+1}}{n+1}\frac{t^n}{n!}\sspeq \frac{1}{t}\, \sum_{n=0}^{\infty}\,C_{n+1}\,\frac{t^{n+1}}{(n+1)!}\sspeq \frac{1}{t}\,({\cal C}(t)\sspm a_0)\, .
\Eeq
We compute the \ogf of the first two terms of the claimed {\sl Faulhaber} formula \Eq{57} multiplied by $d\, (n+1)$
\Beq
d\,(n+1)\,G(d,a;n,x)\sspdef \sum_{m=0}^{\infty}\, x^m\,\left\{ B(d;n+1,x=a+d\,(m+1))\sspm B(d;n+1,x=d)\right\}\,, 
\Eeq
with the polynomials $B(d;n,x)$ from \Eq{52} which are inserted with summation index $k\sspeq n-m$ instead of $m$. The two terms with $k\sspeq n+1$ will be separated and they cancel. Then the sums will be interchanged.
\Beqarray
d\,(n+1)\,G(d,a;n,x)&\sspeq& \sum_{m=0}^{\infty}\, x^m\,\left\{\sum_{k=0}^n\,{\binomial{n+1}{k}}\, B(d;k)\, (a\sspp d\,(m+1))^{n+1-k}\sspm \sum_{k=0}^n\, {\binomial{n+1}{k}}\,B(d;k)\,d^{n+1-k}\right\}\nonumber \\
&\sspeq& \sum_{k=0}^n\,{\binomial{n+1}{k}}\, B(d;k)\,\left\{\left(\sum_{m=0}^{\infty}\,(a\sspp d\,(m+1))^{n+1-k}\,x^m\right) \sspm \frac{1}{1-x}\,d^{n+1-k}\right\}\, .
\Eeqarray
The last term in the curly bracket simplifies with $B(d;k)\sspeq d^k\, B(k)$ from \Eq{53}, and with \Eq{43} rewritten as
\Beq
\sum_{k=0}^n\,{\binomial{n+1}{k}}\,B(k)\sspeq \delta_{n,0}\,,
\Eeq
to  
\dstyle{-\frac{d}{1\sspm x}\,\delta_{n,0}}.\pn
In the remaining double sum one uses the \ogf of sums of powers (see eqs. $(20)$ and  $(21)$) after an index shift $m\sspto m-1$, then one adds and subtracts the new  $m=0$ term. Thus,
\Beqarray
d\,(n+1)\,G(d,a;n,x)&\sspeq& \sum_{k=0}^n\,{\binomial{n+1}{k}}\, B(d;k)\,\frac{1}{x} \left( GP(d,a;n+1-k,x)\sspm a^{n+1-k}\right) \sspm \frac{d}{1\sspm x}\,\delta_{n,0}\nonumber \\
&\sspeq& \sum_{k=0}^n\,{\binomial{n+1}{k}}\, B(d;k)\,\left[-\frac{a^{n+1-k}}{x} \sspp \sum_{m=0}^{n+1-k}\, S2(d,a;n+1-k,m)\,m!\,\frac{x^{m-1}}{(1\sspm x)^{m+1}}\right]\,\nonumber\\
&&\sspm \frac{d}{1\sspm x}\,\delta_{n,0}\, .
\Eeqarray
The term with ${\bf S2}[d,a]$ will now be treated separately as $d\,(n+1)\,G1(d,a;n,x)$ and the remainder\pn
$d\,(n+1)\,G2(d,a;n,x)$ will be added later. In $d\,(n+1)\,G1(d,a;n,x)$ a new summation index $k' \sspeq n+1-k$ is used (called then again $k$), and the $m=0$ sum term will be separated in order to have in both sums the same offset $1$.\psn
\Beqarray
d\,(n+1)\,G1(d,a;n,x) &\sspeq& \sum_{k=1}^{n+1}\,{\binomial{n+1}{k}}\, B(d;n+1-k)\,\left(\sum_{m=1}^k\, S2(d;a;k,m)\,m!\,\frac{x^{m-1}}{(1\sspm x)^{m+1}}\sspp \frac{a^k}{x\,(1\sspm x)}\right)\, \nonumber \\
 &\sspfed&  d\,(n+1)\,G11(d,a;n,x) \sspp  d\,(n+1)\,G12(d,a;n,x)\, .
\Eeqarray
The $m=0$ term $ d\,(n+1)\,G12(d,a;n,x)$ will be added later, and the for the first term we have the following {\it Lemma}.\psn
{\bf Lemma 7:} \psn
\Beq
G11(d,a;n,x)\sspeq GPS(d,a;n,x)\,,
\Eeq
which is the \ogf given in \Eq{2} of the object of desire $PS(d,a;n,m)$, \ie the \lhs of the {\sl Faulhaber} formula \Eq{57}.\psn
{\bf Proof:} The two sums are exchanged, and in the $m$-sum a shift $m\sspto m+1$ will be applied.
\Beqarray
G11(d,a;n,x)&\sspeq& \frac{1}{d\,(n+1)}\,\sum_{m=1}^{n+1}\frac{x^{m-1}}{(1-x)^{m+1}}\,\sum_{k=m}^{n+1}\, {\binomial{n+1}{k}}\,B(d;n+1-k)\,S2(d,a;k,m)\,m!\nonumber\\
&\sspeq&\sum_{m=0}^{n}\frac{x^m}{(1-x)^{m+2}}{\frac{1}{d\,(n+1)}}\sum_{k=m+1}^{n+1} {\binomial{n+1}{k}}B(d;n+1-k)S2(d,a;k,m+1)(m+1)!\,.
\Eeqarray
The $k$-sum will be called $C_{n+1}\sspequiv C(d,a;n+1,m+1)$. Now {\it Lemma 6} is used to compute the \egf of \dstyle{\left\{\frac{C_{n+1}}{d\,(n+1)}\right\}_{n=0}^{\infty}}. Because \dstyle{C_n\sspeq \sum_{k=0}^{n} {\binomial{n}{k}}B(d;n-k)S2(d,a;k,m+1)(m+1)!}  (the sum can start with $k=0$ because $S2(d,a,k,m+1)$ vanishes for $k\sspkl m+1$) is an exponential  convolution, the \egf of $\{C_n\}_{n=0}^{\infty}$ is the product of $EB(d;t)$ from \Eq{55} and the \egf $ES2Col(d,a;t,m+1)$  multiplied by $(m+1)!$, hence 
\Beq
\frac{d\,t}{e^{d\,t}\sspm 1}\cdot e^{a\,t}\,(e^{d\,t}\sspm 1)^{m+1} \sspeq d\,t\,e^{a\,t}\,(e^{d\,t}\sspm 1)^m\, .
\Eeq
Thus, the \egf of \dstyle{\left\{\frac{C_{n+1}}{d\,(n+1)}\right\}_{n=0}^{\infty}} is, by {\it Lemma 6}, $e^{a\,t}\,(e^{d\,t}\sspm 1)^m$ because $C_0 \sspeq C(d,a;0,m+1)\sspeq B(d;0)\,S2(d,a;0,m+1)(m+1)!\sspeq 0$, since $S2(d,a;0,m+1)\sspeq 0$ for $m\sspgeq 0$. But this is the \egf of $\{S2(d,a;n,m)\,m!\}_{n=0}^{\infty}$, and therefore 
\Beq
G11(d,a;n,x) \sspeq \sum_{m=0}^{n}\frac{x^m}{(1-x)^{m+2}}\,S2(d,a;n,m)\,m! \sspeq GPS(d,a;n,x)\, .
\Eeq
In the last step \Eq{8} was used.\hskip 12cm $\square$\psn
If all terms of $G1$ and $G2$ are added we have \psn
\Beqarray
G(d,a;n,x)\sspeq  GPS(d,a;n,x) &\sspp& \frac{1}{d\,(n+1)}\,\left [\frac{1}{x\,(1-x)}\,\sum_{k=1}^{n+1}\,{\binomial{n+1}{k}}\, B(d;n+1-k)\,a^k\right.\nonumber \\
&&\hskip 2cm \left.\sspm \sum_{k=0}^n\,{\binomial{n+1}{k}}\, B(d;k)\,\frac{a^{n+1-k}}{x} \sspm \frac{d}{1-x}\,\delta_{n,0}\right]\ .
\Eeqarray
In the second sum an index change $k'\sspeq n+1-k$ leads to
\Beqarray
G(d,a;n,x) - GPS(d,a;n,x)&=&\frac{1}{d\,(n+1)}\left[\left(\frac{1}{x\,(1-x)}\sspm \frac{1}{x}\right)\sum_{k=1}^{n+1}\,{\binomial{n+1}{k}}B(d;n+1-k)\,a^k\sspm \frac{d}{1-x}\delta_{n,0}\right]\nonumber \\
&\sspeq& \frac{1}{d\,(n+1)}\,\frac{1}{1-x}\,\left[\sum_{k=1}^{n+1}\,{\binomial{n+1}{k}}\, B(d;n+1-k)\,a^k\sspm d\,\delta_{n,0} \right]\nonumber\\
&\sspeq& \frac{1}{d\,(n+1)}\,\frac{1}{1-x}\,\llap{\phantom{xxxxxxx}}\left[\sum_{k=0}^{n+1}\,{\binomial{n+1}{k}}\, B(d;n+1-k)\,a^k\sspm B(d;n+1)\sspm d\,\delta_{n,0} \right]\nonumber \\
&\sspeq& \frac{1}{d\,(n+1)}\,\frac{1}{1-x}\, \left[B(d;n+1,x=a) \sspm B(d;n+1,x=0) \sspm d\,\delta_{n,0}\right].
\Eeqarray
In the last step the polynomials of \Eq{52} have been identified. If now the definition $G(d,a;n,x)$ in \Eq{116} is remembered, and the coefficient $[x^m]$ of this \ogf is picked one finds after this {\it tour de force} the {\it Faulhaber} formula \Eq{57}.\pbn
\psn
{\bf C) Proofs of section 1\,C}\psn
{\bf 1. Proof of eq.\,$\bf (59)$}\psn
In the {\sl Sheffer} group of (infinite) lower triangular matrices the inverse element of ${\bf S}\sspeq (g(x),\,f(x))\sspequiv (g,\,f)$ is 
\Beq
{\bf S}^{-1}\sspeq \left(\frac{1}{g(f^{[-1]}(y))},\,f^{[-1]}(y)\right) \sspequiv \left(\frac{1}{g\circ f^{[-1]}}, \,f^{[-1]}\right)
\Eeq
with the compositional inverse $f^{[-1]}$ of $f$, \ie $f(f^{[-1]}(y))\sspeq y$, or  $f^{[-1]}(f(x))\sspeq x$, identically. Here $f(x)\sspeq x\,{\hat f}(x)$ with ${\hat f}(0)\sspneq 0$. \psn
For the {\sl Sheffer} matrix ${\bf S2}[d,a]$ one has $g(x)\sspeq e^{a\,x}$ and $f(x)\sspeq e^{d\,x}\sspm 1$, hence \dstyle{f^{[-1]}(y)\sspeq \frac{1}{d}\,\log(1\sspp y)}, and 
\Beq
{\bf S1}[d,a]\sspdef ({\bf S2}[d,a])^{-1} \sspeq \left((1\sspp y)^{-\frac{a}{d}} ,\,\frac{1}{d}\,\log(1\sspp y)\right)\,.  
\Eeq
This matrix has in general fractional integer entries. The unsigned matrix ${\bf S1p}[d,a]\sspequiv |{\bf S1}[d,a]|$ ($ p$ for non-negative) has elements $S1p(d,a;n,m)\sspeq (-1)^{n-m}\, S1(d,a;n,m)$ because then the \egf for column $m$ becomes 
\Beq
ES1p(d,a;t,m)\sspeq (1-t)^{-\frac{a}{d}}\,\frac{\left(-\frac{1}{d}\,\log(1\sspm t)\right)^m}{m!}\,,
\Eeq
and both (formal) power series have non-negative elements which are in general fractional numbers.\psn
For combinatorial considerations one is interested in non-negative integer matrices. Therefore, a scaling of the rows is performed: $\widehat{S1p}(d,a;n,m)\sspdef d^n\,S1p(d,a;n,m)$ which leads to diagonal elements $1$, and the {\sl Sheffer} matrix is
\Beq
{\bf \widehat{S1p}}[d,a]\sspeq \left((1\sspm d\,y)^{-\frac{a}{d}},\,-\frac{1}{d}\,\log(1\sspm d\,y)\right)\,,
\Eeq
because the scaling leads to $t\sspto d\,t$ in $ES1p(d,a;t,m)$. The new power series generate exponentially non-negative integers, because \dstyle{\left[\frac{x^n}{n!}\right]\,(1\sspm d\,y)^{-\frac{a}{d}}\sspeq \left(\frac{a}{d}\right)^{\overline{n}}\,d^n \sspeq \prod_{j=0}^{n-1}\,(a\sspp d\,j)\sspeq risefac(d,a;0,n)} (see \Eq{63} for the $risefac$ definition), and  \dstyle{\left[\frac{x^n}{n!}\right]\,\left( -\frac{1}{d}\,\log(1\sspm d\,y)\right)\sspeq (n-1)!\,d^{n-1}} for $n\sspgeq 1$ (and $0$ for $n=0$).\psn 
{\bf 2. Proof of eq.\,$\bf (61)$}\psn
The three term recurrence of the ${\bf \widehat{S1p}}[d,a]$ can be obtained from the \egf of their column sequences $E\widehat{S1p}Col(d,a;t,m)$ given in \Eq{60} which we abbreviate for this proof as $Ep(t,m)$.\psn
{\bf Lemma 8:}\psn
\Beq
(1-d\,t)\,\frac{d\,}{dt}\,Ep(t,m)\sspeq a\,Ep(t,m) \sspp Ep(t,m-1)\, \ \ {\rm for}\ m\sspin \mathbb N\,,
\Eeq
and the input is \dstyle{Ep(t,0)\sspeq (1\sspm d\,t)^{-\frac{a}{d}}}.\psn
{\bf Proof:} This is elementary with \Eq{60}.\psn
Now the recurrence \Eq{61} is seen to satisfy this {\it Lemma}.
\Beqarray
Ep(t,m) &\sspeq& \sum_{n=0}^{\infty}\,\frac{t^n}{n!}\,\widehat{S1p}(d,a;n,m)\nonumber\\
&\sspeq&  \sum_{n=0\ (1)}^{\infty}\,\frac{t^n}{n!}\, \widehat{S1p}(d,a;n-1,m-1)\sspp \sum_{n=0\ (1)}^{\infty}\,\frac{t^n}{n!}\,(d\,n\sspm (d-a))\,\widehat{S1p}(d,a;n-1,m)\nonumber\\
&\sspeq& \sum_{n=0}^{\infty}\,\frac{t^{n+1}}{(n+1)!}\,\widehat{S1p}(d,a;n,m-1) \sspp \sum_{n=0}^{\infty}\,\frac{t^{n+1}}{(n+1)!}\,(d\,(n+1)\sspm (d-a))\,\widehat{S1p}(d,a;n,m)\nonumber \\
&\sspeq& \int dt\, Ep(t,m-1)\sspp d\,t\,Ep(t,m)\sspm (d-a)\, \int dt\,Ep(t,m)\ .
\Eeqarray
In the second line the two sums actually start with $n=1$ because  ${\bf \widehat{S1p}}[d,a]$ vanishes for negative row indices. This is a integral-difference equation with input $Ep(t,0)$ as given in the {\it Lemma}.
\Beq
(1\sspm d\,t)\,Ep(t,m)\sspeq \int dt\,\left(Ep(t,m-1)\sspm (d\sspm a)\,Ep(t,m)\right)\ .
\Eeq
Differentiation produces precisely the equation of the {\it Lemma}.\psn
{\bf 3. Proof of eq.\,$\bf (62)$}\psn
The row polynomials of {\sl Sheffer} triangles are a {\sl Sheffer} transform of the monomials $\{x^n\}_{n=0}^{\infty}$. Therefore, with {\it Lemma 3}, \Eq{82}, the \egf of the (ordinary) row polynomials of ${\bf \widehat{S1p}}[d,a]$ is obtained from the \egf of ${\bf \widehat{S1p}}[d,a]$ given by \Eq{59} and $e^{x\,t}$, \ie
\Beq
EP\widehat{S1p}(d,a;t,x)\sspdef \sum_{n=0}^{\infty}\, P\widehat{S1p}(d,a;n,x)\,\frac{t^n}{n!}\sspeq (1\sspm d\,t)^{-\frac{a}{d}}\,exp(x\,(-\frac{1}{d}\,\log(1-d\,t)))\sspeq (1\sspm d\,t)^{-\frac{x+a}{d}}\, .
\Eeq
Then the binomial theorem leads to 
\Beq
EP\widehat{S1p}(d,a;t,x) \sspeq \sum_{n=0}^{\infty}\,\frac{t^n}{n!}\, (-d)^n\, \left(-\frac{x+a}{d}\right)^{\underline{n}} \sspeq \sum_{n=0}^{\infty}\,\frac{t^n}{n!}\,risefac(d,a;x,n)\, ,
\Eeq
where in the first equation the usual falling factorial $x^{\underline{n}}\sspdef \prod_{j=0}^{n-1}\, (x-j)$ appeared, and in the second equation the definition \Eq{63} for $risefac(d,a;x,n)$ has been used in the rewritten form using ordinary falling factorials.\psn
In this way we have found, as a corollary, the \egf of  $\{risefac(d,a;x,n)\}_{n=0}^{\infty}$ to be \dstyle{(1\sspm d\,t)^{-\frac{x+a}{d}}}.\psn
The $fallfac[d,a]$ analogon is obtained from inverting \Eq{16} using the inverse of the scaled ${\bf \widehat{S2}}[d,a]$  {\sl Sheffer} matrix, \ie the signed ${\bf \widehat{S1}}[d,a]$ matrix.
\Beq
fallfac(d,a;x,n)\sspeq \sum_{m=0}^n\,\widehat{S1}(d,a;n,m)\,x^m\, .
\Eeq
These are the row polynomials of ${\bf \widehat{S1}}[d,a]$.\psn
{\bf 4. Lah[d,a]}\psn
It is tempting to give here the generalized unsigned {\sl Lah} matrix ${\bf L}[d,a]$ as transition matrix between $risefac[d,a]$ and $fallfac[d,a]$.\psn
For the ordinary $[d,a]\sspeq [1,0]$ {\sl Lah} triangle see \seqnum{A271703} (or \seqnum{A008297} with $n\sspgeq m\sspgeq 1$) and \cite{GKP}, exercise 31, p. 312, solution p. 552.\psn
The generalization is
\Beq
 risefac(d,a;x,n)\sspeq \sum_{m=0}^n\,L(d,a;n,m)\,fallfac(d,a;x,m)\ .
\Eeq
From \Eq{62} and \Eq{16} one has, in matrix notation
\Beq
{\bf L}[d,a]\sspeq {\bf\widehat{S1p}}[d,a]\sspcdot {\bf\widehat{S2}}[d,a]\,.
\Eeq
We quote a {\it Lemma} on the multiplication law of the {\sl Sheffer} group.\psn
{\bf Lemma 9:} \psn
If the product of two {\sl Sheffer} matrices with ${\bf S1}\sspeq (g1,\, f1)$ and ${\bf S2}\sspeq (g2,\, f2)$ is ${\bf S3}\sspeq {\bf S1}\sspcdot {\bf S2}$ with ${\bf S3}\sspeq (g3,\, f3)$ then 
\Beq
g3\sspeq g1\, (g2\sspcirc f1)\ \ , \ \ f3\sspeq (f2\sspcirc f1)\ \ ,\ \ {\rm i.e.}, \ \ 
g3(t)\sspeq g1(t)\,g2(f1(t)))\ \ , \ \ f3(t)\sspeq f2(f1(t))\, .
\Eeq
This is standard {\sl Sheffer} lore.\psn
With \Eq{59} and the statement just before \Eq{15} this implies the {\sl Sheffer} structure\psn
\Beq
{\bf L}[d,a]\sspeq \left((1\sspm\,d\,t)^{-\frac{2\,a}{d}},\, \frac{t}{1\sspm d\,t}\right)\, .
\Eeq
The proof of an explicit form along the lines of the mentioned exercise in \cite{GKP} does not immediately lead to an explicit form for $L(d,a;n,m)$ if $a\sspneq 0$. See also the complicated form of \Eq{69} for ${\bf \widehat{S1p}}[d,a]$. Of course the matrix product can  be written with the help of \Eq{64} or \Eq{65} and \Eq{11}.\psn
The \egf of the column sequences is
\Beq
ELCol(d,a;t,m)\sspeq (1\sspm\,d\,t)^{-\frac{2\,a}{d}}\, \frac{1}{m!}\,\left( \frac{t}{1\sspm d\,t} \right)^m\,,\ \ m\sspin \mathbb N_0\, .
\Eeq
From the so called $a-$ and $z-$sequences for {\sl Sheffer} matrices (see the link \cite{WLang2}, where also references are given. This link is found also in \seqnum{A006232}) one finds recurrence relations. The {\it e.g.f.}s of these sequences are ($g$ and $f$ are those of \Eq{140} \psn
\Beqarray
a(y)&\sspeq& \frac{y}{f^{[-1]}(y)} \sspeq 1\sspp d\, y\sspeq a(d;y)\,,\nonumber \\
z(y)&\sspeq& \frac{1}{f^{[-1]}(y)}\,\left(1\sspm \frac{1}{g(f^{[-1]}(y))}\right)\sspeq \frac{1\sspp d\,y}{y}\,\left(1\sspm (1\sspp d\,y)^{-\frac{2\,a}{d}}\right)\sspeq z(d,a;y)\,.
\Eeqarray
This means that there is always a three term recurrence for the matrix entries $L(d,a;n,m)$ for $n \sspgeq m \sspgeq 1$ because the $a-$sequence is $\{1,\,d,\,{\rm repeat}(0)\}$ \ie\psn
\Beq
L(d,a;n,m)\sspeq \frac{n}{m}\,L(d,a;n-1,m-1)\sspp n\, L(d,a;n-1,m)\,\ \ n\sspin \mathbb{N},\, m\sspeq 1,\,2, ...,\, n\, ,
\Eeq
where the input from column $m\sspeq 0$, besides $L(d,a;0,0)\sspeq 1$, can be taken form the \egf \Eq{141}. \psn
In general one can use the $z-$sequence for column $m=0$ in combination with the given recurrence \Eq{143}. For $[d,a]\sspeq [1,0]$ where the $z-$sequence vanishes, the $m \sspeq 0$ column becomes directly $\{1,\,{\rm repeat}(0)\}$. For $[d,a]\sspeq [2,1]$ the $z-$sequence becomes $\{2,\,{\rm repeat}(0)\}$, and the column $m=0$ is also given directly as \seqnum{A000165}. All other cases need also lower row entries with $m\sspgeq 1$.\psn
In this special  ${\bf L}[d,a]$ case one can, however, derive from the column \egf \Eq{114} a four term recurrence, \ie
\Beqarray
L(d,a;n,m)\sspeq L(d,a;n-1,m-1) &\sspp& 2\,(a\sspp d\,(n-1))\, L(d,a;n-1,m)\nonumber \\
&\sspm& d\,(n-1)\,(2\,a\sspp d\,(n-2))\, L(d,a;n-2,m)\,,
\Eeqarray
with inputs $L(d,a;0,0)\sspeq 1$, $L(d,a;n,-1)\sspeq 0$, $L(d,a;-1,m)\sspeq 0$, and $L(d,a;n,m)\sspeq 0$ if $n\sspkl m$.\psn
{\bf Proof:}\,
This uses the definition of \dstyle{ELCol(d,a;t,m)\sspdef \sum_{n=m\, (0)}^{\infty}\, L(d,a;n,m)\,\frac{t^n}{n!}}, and the trivial result 
\Beq
(1-d\,t)^2\,\frac{d\ }{dt}\,ELCol(d,a;t,m)\sspeq 2\,a\,(1\sspm d\, t)\,ELCol(d,a;t,m)\sspp ELCol(d,a;t,m-1)\, .
\Eeq
The recurrence follows then by comparing powers of \dstyle{\frac{t^n}{n!}}\,, sending $n\sspto n-1$.
\pbn
The {\sl Meixner} type recurrence for the row polynomials (see \Eq{85}) is
\Beq
\frac{{\op d}_x}{1\sspp d\,{\bf d}_x}\, PL(d,a;n,x)\sspeq n\, PL(d,a;n-1,x)\,,\ \ n\sspin \mathbb N\,, 
\Eeq
and input $PL(d,a;0,x)\sspeq 1$. The series terminates and this becomes
\Beq
\sum_{k=0}^{n-1}\,(-1)^d\, d^k\,{\bf d}_x^{k+1}\, PL(d,a;n,x)\sspeq n\, PL(d,a;n-1,x)\, .
\Eeq
The general {\sl Sheffer} polynomial recurrence (see \Eq{87} for the rewritten {\sl Roman} corollary) is
\Beq
 PL(d,a;n,x)\sspeq \left((2\,a \sspp x)\,{\bf 1} \sspp 2\,d\,(a\sspp x)\, {\bf d}_x \sspp d^2\,x\,{\bf d}_x^2\right)\, PL(d,a;n-1,x)\,, \ \  n\sspin \mathbb N,\,
\Eeq
and input $PL(d,a;0,x)\sspeq 1$.\pbn
The inverse matrix of ${\bf L}[d,a]$ is also of interest:\psn
\Beq
 fallfac(d,a;x,n)\sspeq \sum_{m=0}^n\,L^{-1}(d,a;n,m)\,risefac(d,a;x,m)\ .
\Eeq
From \Eq{127} one finds the {\sl Sheffer} structure
\Beq
{\bf L}^{-1}[d,a]\sspeq \left(\frac{1}{(1 + d\,t)^{\frac{2\,a}{d}}} ,\, \frac{t}{1\sspp d\,t}\right)\sspeq (gL(-t),\, -fL(-t))\, ,
\Eeq
where $gL$ and $fL$ are taken from \Eq{140}.\psn
This means, by looking at the column {\it e.g.f}.s of {\sl Sheffer} matrices, that the inverse matrix is just obtained by properly signing the ${\bf L}[d,a]$ matrix entries.
\Beq
L^{-1}(d,a;n,\, m)\sspeq (-1)^{n-m}\,L(d,a;n,\,m)\,, \ \ n\sspgeq m\sspgeq 0\,. 
\Eeq
The explicit from of $gL^{-1}[d,a]$ and $fL^{-1}[d,a]$ shows that one has to replace in the ${\bf L}[d,a]$ recurrence formulae $a\sspto -a$ and $d\sspto -d$.\psn
The $a-$ and $z-$sequences are then $aL^{-1}(d)\sspeq \{1,\,-d,\,{\rm repeat}(0)\}$ and the \egf for $zL^{-1}$ is \dstyle{zL^{-1}(d,a;y) \sspeq \frac{1\sspm d\,y}{y}\,\left(1\sspm (1\sspm d\,y)^{-\frac{2\,a}{d}}\right)\sspeq -z(d,a;-y)} (with $z(d,a;y)$ from \Eq{142}). This gives a three term recurrence for $L^{-1}(d,a;n,m)$ for $n\sspgr m\sspgr 1$ with the column $m=0$ as input.\psn
The recurrence derived like above from the column \egf is just \Eq{144} with replacements $a\sspto -a$ and $d\sspto -d$
\Beqarray
L^{-1}(d,a;n,m)\sspeq L^{-1}(d,a;n-1,m-1) &\sspm& 2\,(a\sspp d\,(n-1))\, L^{-1}(d,a;n-1,m)\nonumber \\
&\sspm& d\,(n-1)\,(2\,a\sspp d\,(n-2))\, L^{-1}(d,a;n-2,m)\,,
\Eeqarray
with inputs $L^{-1}(d,a;0,0)\sspeq 1$, $L^{-1}(d,a;n,-1)\sspeq 0$, $L^{-1}(d,a;-1,m)\sspeq 0$, and $L^{-1}(d,a;n,m)\sspeq 0$ if $n\sspkl m$.\psn
The {\sl Meixner} type recurrence for the row polynomials of ${\bf L}^{-1}$ is like the one in \Eq{147} wit replacement $d\sspto -d$
\Beq
\sum_{k=0}^{n-1}\, d^k\,{\bf d}_x^{k+1}\, PL^{-1}(d,a;n,x)\sspeq n\, PL^{-1}(d,a;n-1,x)\, ,
\Eeq
with input $PL^{-1}(d,a;0,x)\sspeq 1$.\psn
The general {\sl Sheffer} recurrence is like \Eq{148} with replacement $a\sspto -a$ and $d \sspto -d$. 
\Beq
 PL^{-1}(d,a;n,x)\sspeq \left((x\sspm 2\,a)\,{\bf 1} \sspm 2\,d\,(x\sspm a)\, {\bf d}_x \sspp d^2\,x\,{\bf d}_x^2\right)\, PL^{-1}(d,a;n-1,x)\,, \ \  n\sspin \mathbb N,\,
\Eeq
and input $PL^{-1}(d,a;0,x)\sspeq 1$.\pbn
{\bf 5. Proof of eq.$\bf(64)$}\psn
It is well known ({\sl Vieta}'s theorem) that the coefficients of a monic polynomial $P(n,\, x)\sspeq \sum_{m=0}^n\, p_m\, x^m$ of degree $n$ are given in terms of the $n$ zeros $x_j,\, j\sspeq 1,\,...,\,n$, of $P$ by $p_m\sspeq (-1)^{n-m}\,\sigma^{(n)}_{n-m}(x_1,x_2,...,x_n)\sspeq -\sigma^{(n)}_{n-m}(-x_1,\,-x_2,\, ...,\,-x_n)$ with the elementary symmetric functions $\sigma^{n}_{n-m}$ of degree $n-m$, and $\sigma^{(n)}_0\sspeq 1$.
For the $risefac(d,a;x,n)$ polynomials \Eq{63} the zeros are $x_j\sspeq -(a\sspp (j-1)\,d) \sspeq -a_{j-1},\, \, j\sspeq 1,\,...,\,n$ proving \Eq{64} for the coefficients $\widehat{S1p}(d,a;n,m)$.\pbn
{\bf 6. Proof of eq. ${\bf (65)}$}\psn
The second version is proved by using \dstyle{ risefac(d,a;x,n)\sspeq d^n\,\left(\frac{x\sspp a}{d}\right)^{\overline{n}}} (see \Eq{63}) and the known result (\eg, \cite{GKP}, p. 263, eq. (6.11)) for ${\bf S1p}$ as transition matrix  
\Beq
x^{\overline{n}}\sspeq \sum_{k=0}^n\,S1p(n,\,k)\,x^k,\ \ k\sspin \mathbb N_0\, .
\Eeq
Now, with the binomial formula,
\Beqarray
risefac(d,a;x,n)&\sspeq& d^n\,\sum_{k=0}^n\,S1p(n,\,k) \left(\frac{x\sspp a}{d}\right)^k\sspeq \sum_{k=0}^n\,d^{n-k}\,S1p(n,\,k)\, \sum_{m=0}^k\,{\binomial{k}{m}}\,a^{k-m}\, x^m \nonumber \\
&\sspeq& \sum_{m=0}^n\,x^m\,\sum_{k=m}^n\,{\binomial{k}{m}}\,S1p(n,\,k)\,a^{k-m}\,d^{n-k}\, ,
\Eeqarray
and the coefficient of $x^m$ is $\widehat{S1p}(d,a;n,m)$ given in the second version of \Eq{65} (with summation index $k\sspto j$). The first version of \Eq{65} is obtained by changing $j\sspto n\sspm j^{\prime}$, using then again $j$ as summation index.\psn
An alternative proof can be given using the recurrence $risfac(d,a;x,n)\sspeq (x\sspp(n-1)\,d\sspp a)\,risefac(d,a;x,n-1)$ with input $risefac(d,a;x,0)\sspeq 1$. Then the {\sl Pascal} recurrence (see \seqnum{A007318}) and the known ${\bf S1p}$ recurrence (given by \Eq{61} for $[d,a]\sspeq [1,0]$) are used.\pbn
{\bf 7. Meixner type recurrence, and proof of eq. $\bf (66)$}\psn
The {\sl Meixner} type recurrence for the monic row polynomials $P\widehat{S1p}(d,a;n,x)$ uses the compositional inverse of the {\sl Sheffer} $f$ function which is \dstyle{\frac{1}{d}\,(1\sspm e^{-d\,y})} (see \Eq{85}). Therefore, $f^{[-1]}({\bf d}_x)\, P\widehat{S1p}(d,a;n,x)\sspeq n\, P\widehat{S1p}(d,a;n-1,x)$ becomes
\Beq
\sum_{k=1}^n\, (-1)^{k-1}\,\frac{d^{k-1}}{k!}\,({\bf d}_x)^k\, P\widehat{S1p}(d,a;n,x)\sspeq n\, P\widehat{S1p}(d,a;n-1,x)\,, 
\Eeq
with input $P\widehat{S1p}(d,a;0,x)\sspeq 1$.\psn
The general {\sl Sheffer} recurrence (see {\sl Lemma} $5$, \Eq{87}, also for the {\sl Roman} reference) uses, with the {\sl Sheffer} $g$ and the given $f^{[-1]}$ function,
\Beq
g(f^{[-1]}(t))\sspeq e^{a\, t},\ \ \frac{d}{dt}\,g(f^{[-1]}(t))\sspeq a\, e^{a\,t},\ \  \frac{d}{dt}\,f^{[-1]}(t)\sspeq e^{-d\,t}\, ,
\Eeq
leading to 
\Beq
P\widehat{S1p}(d,a;n,x)\sspeq (x\sspp a)\,e^{d\,{\bf d_x}}\,P\widehat{S1p}(d,a;n-1,x)\sspeq (x\sspp a)\, P\widehat{S1p}(d,a;n-1,x+d),\,
\Eeq
by {\sl Taylor}'s theorem, and this proves \Eq{66}.\pbn
{\bf Proof of eq, ${\bf (67)}$}\psn
This is covered by {\it Corollary 1} (after \Eq{82}).\pbn
{\bf Proof of eq. ${\bf(69)}$}\psn
The direct way uses the inversion of \Eq{18}.\psn
{\bf Lemma 10:}\psn
\Beq
S2(n,\,m) \sspeq \left(\frac{-a}{d}\right)^n\,\sum_{k=0}^n\,(-1)^k\,{\binomial{n}{k}}\,a^{-k}\,S2(d,a;k,m)\,, \ \ {\rm for}\ \ n\sspgeq m\sspgeq 0\ .
\Eeq
{\bf Proof}: From the exponential convolution \Eq{18} one has for the column \egf \dstyle{ES2Col(d,a;t,m)\sspeq e^{a\,t}\, ES2Col(d\,t,m)} (see \Eq{13}). This means (for $d\sspin \mathbb N $) that
\Beq
ES2Col(t,m)\sspeq e^{-\frac{a}{d}\,t}\, ES2Col\left(d,a;\frac{t}{d},m\right)\, ,
\Eeq
which gives for $S2(n,\,m)$ the exponential convolution leading to the assertion.\hskip 5.5cm $\square$
\psn
Then the proof of \Eq{69} starts by inserting $S1p(n\, j)$, in the second version of \Eq{65}, by {\sl Schl\"ohmilch}'s formula, \Eq{68}, with $S2(n-j+k,\,k)$ from the {\it Lemma}, \Eq{160}:
\Beqarray
\widehat{S1p}(d,a;n,m)&\sspeq& \sum_{j=m}^n\,{\binomial{j}{m}}\,a^{j-m}\,d^{n-j}\,(-1)^{n-j}\,\sum_{k=0}^{n-j}\,(-1)^k\,{\binomial{n+k-1}{j-1}}\,{\binomial{2\,n-j}{n-j-k}}\,* \nonumber \\
&& *\, \left(\frac{-a}{d}\right)^{n-j+k}\ \sum_{l=0}^{n-j+k}\,(-1)^l\,{\binomial{n+k-j}{l}}\,a^{-l}\,d^k\,\widehat{S2}(d,a;l,k)\, ,
\Eeqarray
where ${\bf S2}[d,a]$ has been replaced by ${\bf \widehat{S2}}[d,a]$ (see the line before \Eq{15}). Collecting $a$ and $d$ powers and the signs leads to \Eq{69}.\psn
Another more complicated proof follows the one of the usual {\sl Schl\"ohmilch} formula given in \cite{Charalambides}, p. 290. This uses the {\sl Lagrange} inversion theorem for powers of a (formal) power series.\psn
{\bf Lemma 11}: {\bf Lagrange theorem and inversion} \cite{Fichtenholz}, p. 523, \Eq{29},  \cite{WhittakerWatson}, p. 133.\psn
{\bf a)} With ${\tilde f}(x) = f(y(x))$, $y(x)\sspeq a\sspp x\,\varphi(y)$  (here as formal power series)
\Beq
{\tilde f}(x)\sspeq f(a)\sspp \sum_{n=1}^{\infty}\,\frac{x^n}{n!} \frac{d^{n-1}\ }{da^{n-1}}\left[\varphi^n(a)\,f^{\prime}(a)\right]\, .
\Eeq
{\bf b)} With $a\sspeq 0$, $y\sspeq x\,\psi(x)$, and $f(y)\sspeq y^k$, $k\sspin \mathbb N_0$ and $x(y)\sspeq y^{[-1]}(y)$ (compositional inverse)\psn
\Beqarray
x^k(y)&\sspeq& \delta_{k,0}\sspp k\,\sum_{n=1}^{\infty}\, \frac{y^n}{n!}\,\left.\frac{d^{n-1}\ }{dt^{n-1}}\left[\left(\frac{1}{\psi(t)}\right)^n\,t^{k-1}\right]\right|_{t=0}\nonumber\\
&\sspeq& \delta_{k,0} \sspp \sum_{n=1}^{\infty}\,\frac{y^n}{n!}\,\sum_{j=0}^{n-1}\, {\binomial{n-1}{j}}\,k^{\underline{n-j}}\,\left.\left[\frac{d^j\ }{dt^j}\left(\frac{1}{\psi(t)}\right)^n\right]\,t^{k-n+j}\right|_{t=0}\nonumber \\
&\sspeq& \delta_{k,0} \sspp k!\,\sum_{m=k}^{\infty}\,\frac{y^m}{m!}\,{\binomial{m-1}{m-k}}\,\left.\left[\frac{d^{m-k}\ }{dt^{m-k}}\left(\frac{1}{\psi(t)}\right)^m\right]\right|_{t=0}\, .
\Eeqarray
{\bf Proof}: Part {\bf a)} is the standard {\sl Lagrange} theorem.\psn
The first equation of part {\bf b)} follows by exchanging the r\^ole of $y$ and $x$, using \dstyle{\varphi(x)\sspeq \frac{1}{\psi(x)}} and $0^k\sspeq \delta_{k,0}$. (See \cite{Fichtenholz}, pp. 524-525 for the case $k=1$). The second equation uses the {\sl Leibniz} rule. Then only $j\sspeq n-k\sspgeq 0$ survives after evaluation at $t=0$, and in the last formula the summation index has been changed for later purposes from $n$ to $m$.\pbn
{\bf Corollary 2}:\psn
\Beq
\left.\frac{d^n\ }{dy^n}\left[\frac{1}{k!}\,(y^{[-1]}(y))^k\right]\right|_{y=0}\sspeq {\binomial{n-1}{n-k}}\, \left.\frac{d^{n-k}\ }{dt^{n-k}}\left[\frac{1}{\psi^n(t)}\right]\right|_{t=0}\,, \ \ {\rm for}\ \ n\sspgeq k\sspin \mathbb N  .
\Eeq
The $\delta$ term now disappeared for $k\sspgeq 1$.\pn
This coincides with the inversion formula of {\sl Lagrange} given in \cite{Charalambides}, Theorem 11.11, p. 435, used in the proof on p. 290. \pbn
Another formula is needed to convert later the negative powers of $\psi$ in the {\sl Corollary} into positive ones. This is given in \cite{Charalambides}, as {\sl Remark 11.5}, p. 432, also used in the proof on p. 290.
\Beq
\left.\frac{d^m\ }{dt^m}\left[(h(t))^s\right]\right|_{t=0}\sspeq \sum_{r=0}^m\,{\binomial{s}{r}}\,{\binomial{m-s}{m-r}}\,\left.\frac{d^m\ }{dt^m}\left[(h(t))^r\right]\right|_{t=0}\,,\ \ {\rm for}\ h(0)\sspeq 1\,, \ m \sspin \mathbb N_0,\ s\sspin \mathbb R\,.   
\Eeq
Now we start with the derivation of a generalized {\sl Schl\"ohmilch} formula (we use here the column index $k$).
\Beq
\widehat{S1p}(d,a;n,k)\sspeq \left.\frac{d^n\ }{dy^n}\left[ (1\sspm d\, y)^{-\frac{a}{d}}\,\frac{1}{k!}\,\left(-\frac{1}{d}\,\log(1\sspm d\,y)\right)^k \right]\right|_{y=0}\,\  {\rm for}\ \ 0\sspleq k\sspleq n\, .
\Eeq
The result for $k\sspeq 0$ is known from \Eq{62} to be $risefac(d,a;x=0,n)$. The {\sl Leibniz} rule is applied (because there is no closed formula for the inverse of the product in the solid brackets, written as $k-$th power). The derivatives of the first factor are known (see the $k\sspeq 0$ result) as
\Beq
\left.\frac{d^{n-m}\ }{dy^{n-m}}\left[ (1\sspm d\, y)^{-\frac{a}{d}}\right]\right|_{y=0}\sspeq risefac(d,a;x=0,n-m)\, . 
\Eeq
For the second factor the {\sl Corollary} is applied with $n\sspto m$, and the known compositional inverse of \dstyle{y^{[-1]}(d;y)\sspeq -\frac{1}{d}\,\log(1\sspm d\,y)} \viz \dstyle{y(d;t)\sspeq t\,\psi(d;t)\sspeq \frac{1}{d}\,(1\sspm e^{-d\,t})}, is used.
\Beq
\left. \frac{d^{m}\ }{dy^{m}}\left[\frac{1}{k!}\,\left(-\frac{1}{d}\,\log(1\sspm d\,y)\right)^k\right]\right|_{y=0} \sspeq {\binomial{m-1}{m-k}}\,\left. \frac{d^{m-k}\ }{dt^{m-k}}\left[\left(\frac{1\sspm e^{-d\,t}}{d\,t}\right)^{-m}\right]\right|_{t=0}\, .
\Eeq
The negative power on the \rhs is now converted with the help of \Eq{166} with $s\sspto -m$ and $m\sspto m-k$. Then the binomial with negative upper entry is rewritten as \dstyle{{\binomial{-m}{r}}\sspeq (-1)^r\,{\binomial{r+m-1}{r}}} (see \Eq{95}), and we use the abbreviation \dstyle{\psi(d;t)\sspeq \frac{1\sspm e^{-d\,t}}{d\,t}} from above.
\Beq
\left. \frac{d^{m}\ }{dy^{m}}\left[\frac{1}{k!}\,\left(-\frac{1}{d}\,\log(1\sspm d\,y)\right)^k\right]\right|_{y=0} \sspeq {\binomial{m-1}{m-k}}\,\sum_{r=0}^{m-k}\,(-1)^r\,{\binomial{r+m-1}{r}}\,{\binomial{2\,m-k}{m-k-r}}\, \left. \frac{d^{m-k}\ }{dt^{m-k}}\left(\psi(d;t)^r\right)\right|_{t=0}\, . 
\Eeq
With the \egf \dstyle{E{\widehat{S2}}Col(d,a;x,r)\sspeq e^{a\,x}\,\frac{\left( \frac{1}{d}\,(e^{d\,x}\sspm 1)\right)^r}{r!}} for column $r$ of $\bf{\widehat{S2}}[d,a]$ from its {\sl Sheffer} structure one obtains after multiplication with $x^{-r}$ and a sign flip $x\sspeq -t$ 
\Beq
\psi(d;t)^r\sspeq \left(\frac{1\sspm e^{-d\,t}}{d\,t}\right)^r\sspeq e^{a\,t}\,\sum_{n=r}^{\infty}\,(-1)^{n-r}\,\widehat{S2}(d,a;n,r)\,r!\frac{t^{n-r}}{n!}\sspfed e^{a\,t}\,A(d,a;t,r).
\Eeq
For the $(m-k)$-th derivative evaluated at $t\sspeq 0$ the {\sl Leibniz} rule is again applied with a summation index $p$. The exponential factor leads to $a^{m-k-p}$. The $p$-th derivative w.r.t. $t$ of $A(d,a;t,r)$ is evaluated at $t=0$ after the index shift $n-r\sspeq s$ in $A$. This leads to the collapse of the $s-$sum due to the $t\sspeq 0$ evaluation, whence $s\sspeq p$, and $p^{\underline p}\sspeq p!$. The result is
\Beq
\left.\frac{d^p\ }{dt^p}\, A(d,a;t,r) \right|_{t=0}\speq \frac{r!\,p!}{(p+r)!}\, (-1)^p\,\widehat{S2}(d,a;p+r,r)\, .
\Eeq
Thus
\Beq
\left. \frac{d^{m-k}\ }{dt^{m-k}}\left(\psi(d;t)^r\right)\right|_{t=0} \sspeq \sum_{p=0}^{m-k}\,{\binomial{m-k}{p}}\,a^{m-k-p}\,(-1)^p\,\frac{1}{{\binomial{p+r}{p}}}\,\widehat{S2}(d,a,p+r,r) \,,\ \ {\rm for}\ \ m\sspgeq k\,.
\Eeq
This leads, with the abbreviation $risefac(d,a,x=0,n-m) = (d,a)^{\overline{n-m}}$, to
\Beqarray
\widehat{S1p}(d,a;n,k)&\sspeq& \sum_{m=k}^n\,{\binomial{n}{m}}\,(d,a)^{\overline{n-m}}\,{\binomial{m-1}{m-k}}\sum_{r=0}^{m-k}\, (-1)^r\, {\binomial{r+m-1}{r}}\,{\binomial{2\,m-k}{m+r}}\, * \nonumber \\
&& * \sum_{p=0}^{m-k}\,(-1)^p\,\frac{{\binomial{m-k}{p}}}{{\binomial{p+r}{r}}}\, a^{m-k-p}\, \widehat{S2}(d,a;p+r,r)\,,\ \ {\rm for}\  n\sspgeq m\sspgeq 1,\,   
\Eeqarray
A new summation index $m-k\sspeq m'$ is used, and the sum over the triangular array, with rows indexed by $m$ and columns by $r$, is reordered by summing first the columns $r\sspeq 0,\,...,\,n-k$ and then the rows $m\sspeq r,\,...,\,n-k$. The binomials in these two sums are rewritten and the final result is (using column index $m$ instead of $k$)
\Beqarray
{\fbox{\color{bleudefrance}$\widehat{S1p}(d,a;n,m)$}}&\sspeq& \frac{n!}{(m-1)!} \sum_{r=0}^{n-m}\,\frac{(-1)^r}{r!}\, \sum_{k=r}^{n-m}\, (d,a)^{\overline{n-k-m}}\, {\binomial{2\,k+m}{k+m}}\, \frac{1}{k+m+r}\,\frac{1}{(n-m-k)!}\, \frac{1}{(k-r)!}\, * \nonumber \\
&& * \sum_{p=0}^k\,(-1)^p\,\frac{{\binomial{k}{p}}}{{\binomial{p+r}{r}}}\, a^{k-p}\, \widehat{S2}(d,a;p+r,r)\,\ {\rm for}\  n\sspgeq m\sspgeq 1,\,   
\Eeqarray
and $\widehat{S1p}(d,a;n,0)\sspeq risefac(d,a;x=0,n)$.
\pbn
The binomials could be rewritten by building some binomials, but no essential simplification seems to be possible.\psn
Both generalizations of the {\sl Schl\"ohmilch} formula involve three sums, but the direct version given in \Eq{69} looks simpler than the second version \Eq{175}.
\pbn
The statements in {\it section D} are obvious.
\pbn 

\pbn
\hrulefill
\pbn
{\it 2010 Mathematics Subject Classification}: Primary 11B68, 11B73, Secondary  	05A15, 11B25, 11B37.\psn
{Keywords}: Sums of powers of positive integers, arithmetic progression, generating functions, Stirling numbers, Lah numbers, Eulerian numbers, Bernoulli numbers.
\pbn
\hrulefill
\pbn
{\it OEIS} \cite{OEIS} A numbers: \psn
\seqnum{A000142}, \seqnum{A000165}, \seqnum{A001147}, \seqnum{A006232}, 
\seqnum{A007559}, \seqnum{A007318}, \seqnum{A007559}, \seqnum{A007696}, \seqnum{A008297}, \seqnum{A008544}, \seqnum{A008545}, \seqnum{A008546}, \seqnum{A008548}, \seqnum{A027641}, \seqnum{A027642}, \seqnum{A028338}, \seqnum{A032031}, \seqnum{A039755}, \seqnum{A047053}, \seqnum{A048993}, \seqnum{A052562}, \seqnum{A053382}, \seqnum{A053383}, \seqnum{A060187}, \seqnum{A123125}, \seqnum{A131689}, \seqnum{A132393}, \seqnum{A141459}, \seqnum{A143395}, \seqnum{A145901}, \seqnum{A154537}, \seqnum{A157779}, \seqnum{A157799}, \seqnum{A157817}, \seqnum{A157833}, \seqnum{A157866}, \seqnum{A173018}, \seqnum{A196838}, \seqnum{A196839}, \seqnum{A225117}, \seqnum{A225118}, \seqnum{A225466}, \seqnum{A225467}, \seqnum{A225468}, \seqnum{A225472}, \seqnum{A225473}, \seqnum{A239275}, \seqnum{A271703}, \seqnum{A282629}, \seqnum{A284861}, \seqnum{A285061}, \seqnum{A285066}, \seqnum{A285068}, \seqnum{A285863}, \seqnum{A286718}, \seqnum{A288872}, \seqnum{A288873}.
\pbn
\hrulefill
\pbn

\begin{thebibliography}{99}
\bibitem{Bala} P. Bala, A 3 parameter family of generalized Stirling numbers, \url{https://oeis.org/A143395/a143395.pdf}.This link appears, e.g., under \seqnum{A143395}.
\bibitem{Brenke} W. C. Brenke, On generating functions of polynomial systems, Am. Math. Monthly 52 (1945) 297-301.
\bibitem{Charalambides}  Ch. A. Charalambides, {\em Enumerative Combinatorics}, Chapman \& Hall/CRC, 2002.
\bibitem{Chihara} T. S. Chihara, {\em An Introduction to Orthogonal Polynomials}, Gordon and Breach, New York, London, Paris, 1978.
\bibitem{CoGuy} Johm H. Conway and Richard K, Guy, {\em The Book of Numbers}, Copernicus, an Imprint of Springer-Verlag, 1995.
\bibitem{Edwards1} A. W. F. Edwards, Sums of powers of integers: a little of the history, Mat. Gazette Vol. 66, No. 435 (1982) 22 - 28.
\bibitem{Edwards2} A. W. F. Edwards, A quick route to sums of powers, Am. Math. Monthly, Vol. 93, No. 6 (1986) 451 - 455.
\bibitem{Fichtenholz} G. M. Fichtenholz, {\em Differentail- und Integralrechnung} II, p. 523, VEB Deutscher Verlag der Wissenschaften, Berlin, 1964.
\bibitem{GKP} R. L. Graham, D. E. Knuth and O. Patashnik, {\em Concrete Mathematics}, 2nd edition, 1994, Addison-Weseley, Reading, Massachusetts, 1991.
\bibitem{Hardy} G.H. Hardy, {\em Ramanujan, Twelve lectures...}, AMS Chelsea Publishing, Providence, Rhode Island, 2002. 
\bibitem{Hawlitschek} Kurt Hawlitschek, {\em Jahann Faulhaber 1580 - 1635, Eine Bl\"utezeit der mathematischen Wissenschaften in Ulm}, Ver\"offentlichungen der Universit\"at Ulm, Band 18, Stadtbibliothek Ulm 1995.
\bibitem{Knuth} Donald E. Knuth, Johann Faulhaber and Sums of Powers, Math. of Comp., Vol. 61, No. 203 (1993) 277 - 194.
\bibitem{Koecher} Max Koecher, {\em Klassische elementare Analysis}, Birkh\"auser Verlag, Basel $\cdot$ Boston, 1987.
\bibitem{Krishnamurthy} V. Krishnamurthy, {\em Combinatorics, theory and application}, Ellis Horwood Limited, 1986. 
\bibitem{WLang} Wolfdieter Lang, Combinatorial Interpretation of Generalized Stirling Numbers, Journal of Integer Sequences, Vol. 12 (2009), Article 09.3.3, \url{https://cs.uwaterloo.ca/journals/JIS/VOL12/Lang/lang.html}.
\bibitem{WLang2} Wolfdieter Lang, Sheffer $a-$ and $z-$sequences, \url{https://oeis.org/A006232/a006232.pdf}. 
\bibitem{Luschny1} Peter Luschny, The Stirling-Frobenius numbers, \pn
\url{http://www.luschny.de/math/euler/StirlingFrobeniusNumbers.html}. This link appears, e.g., under \seqnum{A225466}.
\bibitem{Luschny2} Peter Luschny, Eulerian polynomials,\pn \url{http://www.luschny.de/math/euler/GeneralizedEulerianPolynomials.html#generalization_of_the_Eulerian_polynomials}. 
\pn This link appears, e.g., under \seqnum{A225466}.
\bibitem{Maple}  Maple$^{TM}$, \url{http://www.maplesoft.com/}.  
\bibitem{Meixner} J. Meixner, Orthogonale Polynomsysteme mit einer besonderen Gestalt der erzeugenden Funktion, J. London Mat. Soc. 9 (19334) 6-13. 
\bibitem{OEIS} The On-Line Encyclopedia of Integer Sequences (2010), published electronically at \url{http://oeis.org}.
\bibitem{Roman} Steven Roman, {\em The Umbral Calculus}, Academic Press, London, 1984.
\bibitem{Schneider} Ivo Schneider, {\em  Johannes Faulhaber,1580 - 1635, Rechenmeister in einer Welt des Umbruchs}, Birkh\"auser, Basel $\cdot$ Boston $\cdot$ Berlin, 1993.
\bibitem{WhittakerWatson} E. T. Whittaker and G. N. Watson, {\em A Course of Modern Analysis}, Fourth ed., Cambridge, at the University Press, 1958, p. 133. 
\end{thebibliography}
\end{document}